\documentclass[journal]{IEEEtran}

\usepackage{etex}
\reserveinserts{28}

\usepackage{cite}
\usepackage{amsmath}
\usepackage{amssymb}
\usepackage{graphicx}
\usepackage{subfigure}
\usepackage{epstopdf}
\usepackage{amsfonts}
\usepackage[all]{xy}
\usepackage{array,xcolor}
\usepackage{color}

\usepackage{caption}
\usepackage{algorithm}
\usepackage{algpseudocode}
\usepackage{multirow}

\usepackage{threeparttable}
\usepackage{booktabs}
\usepackage{changepage}

\usepackage{stfloats}
\usepackage{balance}

\usepackage{tikz}
\usetikzlibrary{shapes,matrix,decorations.shapes}
\usetikzlibrary{arrows,decorations.pathmorphing,backgrounds,positioning,fit,matrix}
\usetikzlibrary{shapes,arrows,chains}
\usetikzlibrary[calc]
\usetikzlibrary[plotmarks]
\usetikzlibrary[patterns,decorations]
\tikzset{fontscale/.style = {font=\relsize{#1}}
    }

\newtheorem{myDef}{Definition}
\newtheorem{myTheo}{Theorem}

\newtheorem{myRem}{Remark}
\newtheorem{myLem}{Lemma}

\newcounter{tempEquationCounter}
\newcounter{thisEquationNumber}
\newenvironment{floatEq}
{\setcounter{thisEquationNumber}{\value{equation}}\addtocounter{equation}{1}
\begin{figure*}[t]
\normalsize\setcounter{tempEquationCounter}{\value{equation}}
\setcounter{equation}{\value{thisEquationNumber}}
}
{\setcounter{equation}{\value{tempEquationCounter}}
\hrulefill\vspace*{4pt}
\end{figure*}
}

\begin{document}
\title{Cooperative Control of Heterogeneous Connected Vehicles with Directed Acyclic Interactions}

\author{Yang~Zheng,~\IEEEmembership{Student Member,~IEEE}, Yougang~Bian, Shen Li, \\ and Shengbo~Eben~Li,~\IEEEmembership{Senior Member,~IEEE}
\thanks{This research is supported by NSF China with 51575293 and 51622504, National Key R\&D Program of China with 2016YFB0100906, and International Sci\&Tech Cooperation Program of China under 2016YFE0102200. All correspondence should be sent to S. E. Li.}
\thanks{Y. Zheng was with the Department of Automotive Engineering at Tsinghua University. He is now with the Department of Engineering Science at the University of Oxford. (e-mail: yang.zheng@eng.ox.ac.uk).}
\thanks{Y. Bian and S. E. Li are with the State Key Lab of Automotive Safety and Energy, Department of Automotive Engineering, Tsinghua University, Beijing, 100084, China. (e-mail: byg14@mails.tsinghua.edu.cn, lishbo@tsinghua.edu.cn).}
\thanks{S. Li is with the Department of Civil \& Environmental Engineering at University of Wisconsin-Madison (e-mail: shen.li@wisc.edu).}
}

\maketitle

\begin{abstract}
  Cooperation of multiple connected vehicles has the potential to benefit the road traffic greatly. In this paper, we consider analysis and synthesis problems of the cooperative control of a platoon of heterogeneous connected vehicles with directed acyclic interactions (characterized by directed acyclic graphs). In contrast to previous works that view heterogeneity as a type of uncertainty, this paper directly takes heterogeneity into account in the problem formulation, allowing us to develop a deeper understanding of the influence of heterogeneity on the collective behavior of a platoon of connected vehicles. Our major strategies include an application of the celebrated internal model principle and an exploitation of lower-triangular structures for platoons with directed acyclic interactions. The major findings include: 1) we explicitly highlight the tracking ability of heterogeneous platoons, showing that the followers can only track the leader's spacing and velocity; 2) we analytically derive a stability region of feedback gains for platoons with directed acyclic interactions; 3) and consequently we propose a synthesis method based on the solution to an algebraic Riccati equation that shares the dimension of single vehicle dynamics. Numerical experiments are carried out to validate the effectiveness of our results.
\end{abstract}


\IEEEpeerreviewmaketitle

\section{Introduction}\label{section:1}
\IEEEPARstart{C}{onnected} vehicles have recently received increasing attention from both academia and industry, due to the high potential to significantly benefit road transportation~\cite{fagnant2015preparing,li2017dynamical,englund2016grand, zheng2015dynamic}. One important application is to develop \emph{cooperative control strategies} for multiple connected vehicles based on local information to guarantee a certain global coordination, such that the traffic efficiency and road safety are improved. This is called \emph{platoon control} of connected vehicles~\cite{li2017dynamical}, and is also known as \emph{cooperative adaptive cruise control} (CACC)~\cite{oncu2014cooperative}.

One main objective in the cooperative control of multiple connected vehicles is to ensure that all vehicles in a group maintain a desired cruising velocity while keeping a pre-specified inter-vehicle distance. This problem has a long history in control theory dating back to the pioneering work in the 1960s~\cite{levine1966optimal}, where an optimal control framework was introduced to deal with spacing regulation of multiple vehicles. The earliest practices on platoon control can be traced back to the California Partners for Advanced Transportation Technology (PATH) program in the 1980s~\cite{shladover1991automated}, where many practical issues were discussed, including control architecture, spacing policies, sensors and actuators, and string stability. Recent advances have emerged in the application of advanced control methods for platooning of connected vehicles, such as $\mathcal{H}_{\infty}$ control~\cite{ploeg2014controller}, distributed model predictive control~\cite{dunbar2012distributed,Zheng2016distributed}, and sliding mode control~\cite{guo2016distributed}. Also, some proof-of-concept demonstrations have been performed in the projects of GCDC~\cite{englund2016grand}, SARTRE~\cite{robinson2010operating} and Energy-ITS~\cite{tsugawa2011automated}. The interested reader is referred to~\cite{li2017dynamical} for a recent overview.

With the rapid development of vehicle-to-vehicle (V2V) communications, such as DSRC and VANETs~\cite{willke2009survey}, one recent research focus of platooning is on developing scalable analysis and synthesis methods for the cooperation of large-scale platoons with various communication topologies~\cite{li2017dynamical}. For instance, an explicit stabilizing region of linear feedback gains was derived for homogenous platoons with a large class of communication topologies~\cite{zheng2016stability}. Barooah \emph{et al.} introduced a mistuning-based design method to improve closed-loop stability margin of platoons with bidirectional topologies~\cite{barooah2009mistuning}, which has recently been extended to cover the inertial time lag of vehicle powertrains in~\cite{zheng2015stabilityMargin}. One tradeoff of the mistuning-based controller was highlighted in~\cite{herman2015nonzero}, where it is shown that the closed-loop $\mathcal{H}_{\infty}$ norm increases exponentially as the platoon size grows. This fact is consistent with a high-gain condition in the design of distributed $\mathcal{H}_{\infty}$ controllers for platoons with undirected topologies~\cite{zheng2017platooning}. More recently, Qin and Orosz proposed a decomposition method for scalable stability analysis of large connected vehicle systems, where stochastic communication delays were covered~\cite{qin2017scalable}. The aforementioned studies have offered efficient methods for performance analysis and controller synthesis of large-scale platoons of connected vehicles. However, most of them require a key assumption that the dynamics of each vehicle are homogeneous~\cite{zheng2016stability,barooah2009mistuning,zheng2015stabilityMargin,herman2015nonzero,zheng2017platooning,qin2017scalable}. This assumption allows one to apply the decomposition results in multi-agent systems~\cite{fax2004information,li2010consensus} and greatly simplifies the theoretical analysis and synthesis.

The assumption of homogeneity may be too restrictive and impractical since diverse types of vehicles should be allowed in a platoon formation. This leads to the design of heterogeneous platoons, which actually attracts research attention as early as the practices in the PATH program~\cite{shladover1991automated}. For instance, an inclusion principle was applied to decompose a string of interconnected heterogeneous vehicles into a set of subsystems, in which overlapping controllers were designed~\cite{stankovic2000decentralized}. String stability analysis for heterogeneous vehicle platoons was discussed in~\cite{shaw2007string} and~\cite{lestas2007scalability}. Naus \emph{et al.} derived a necessary and sufficient condition in frequency-domain for string stability, where heterogeneous traffic was taken into account~\cite{naus2010string}. In~\cite{ploeg2014lp}, Ploeg \emph{et al.} introduced an $\mathcal{L}_{p}$ string stability of cascade systems using input-output properties, which is suitable for heterogeneous platoons with nearest interactions. In~\cite{harfouch2017adaptive}, an adaptive switched control approach was proposed for heterogeneous platoons with communication losses. Besides, R{\"o}d{\"o}nyi discussed an adaptive spacing policy that is able to guarantee string stability~\cite{rodonyi2017adaptive}. The results of~\cite{shaw2007string, lestas2007scalability, naus2010string,ploeg2014lp,harfouch2017adaptive,rodonyi2017adaptive} offer some insights on the design of heterogeneous platoons. However, these methods are only applicable for very limited types of communication topologies, \emph{e.g.}, predecessor-following (PF) type and predecessor-leader following (PLF) type, since most of them rely on the exploitation of a cascade structure in an implicit or explicit way. Recently, many other types of topologies are emerging, \emph{e.g.}, the multiple-PF type~\cite{zheng2016stability}, thanks to the rapid deployment of V2V techniques. New challenges naturally arise for cooperative control of heterogeneous connected vehicles considering the variety of topologies. Note that there are a few recent works that try to address this issue; see \emph{e.g.},~\cite{gao2016robust,li2012distributed,kim2011output,baldi2018adaptive}. In these works, the authors typically consider the heterogeneity as a type of uncertainty and assume that the vehicles share an identical nominal model. It means that the results on homogeneous platoons~\cite{zheng2017platooning,qin2017scalable,fax2004information,li2010consensus} can be basically applied to the nominal platoon system.

In this paper, we consider cooperative control of a platoon of heterogeneous connected vehicles with directed acyclic interactions (see the precise definition in Section~\ref{section:4b}), and directly take the heterogeneity into account. In contrast to previous studies, this treatment allows us to develop a deeper understanding of the effects of heterogeneity on the collective behavior of a platoon, as well as to highlight the influence of topological variety introduced by V2V techniques. Under the notion of directed acyclic interactions, the closed-loop heterogeneous platoon system becomes decomposable thanks to a lower-triangular structure. This technique does not rely on the normal eigenvalue decomposition or similarity transformation that is widely used in homogeneous platoons (see~\cite{li2017dynamical,li2015overview} for example). More precisely, our contributions are:
\vspace{-0.5mm}
\begin{enumerate}
  \item The internal model principle, a fundamental result in heterogeneous multi-agent systems~\cite{wieland2011internal}, is applied to the analysis of heterogeneous platoons. The internal model principle presents a necessary and sufficient condition for the synchronization of heterogeneous linear networks. We discuss the implication of internal model principle for cooperative control of heterogeneous connected vehicles, and highlight the tracking ability of each following vehicle. To reach the steady consensus state, the leader should run at a constant speed, indicating that the followers can only track the leader's spacing and velocity.
  \item We derive an explicit and analytical region for feedback gains that guarantee the asymptotical stability of heterogeneous platoons with directed acyclic interactions. This result not only explicitly highlights the necessity of spacing and velocity information to stabilize a heterogeneous platoon, but also points out that the existence of a spanning tree in the communication graph is essential for stabilization. The influence of the heterogeneity in vehicle dynamics is directly reflected in the stability region. Our result generalizes the stability condition in~\cite{zheng2016stability} to heterogeneous platoons.
  \item According to the internal stability result, we propose a synthesis method based on the solution to an algebraic Riccati equation (ARE) relying on the dynamics of each individual vehicle. By exploiting a lower-triangular structure, this design method is decoupled from the communication graph. This makes the computational complexity independent of the platoon size. Besides, the synthesis method has a relatively clear physical interpretation on the convergence rate design, which  facilitates its application in practice.
\end{enumerate}

The rest of this paper is organized as follows. Section~\ref{section:2} presents the problem statement. The tracking ability of a heterogeneous platoon is discussed using the internal model principle in Section~\ref{section:3}. Section~\ref{section:4} presents the stability result and introduces an ARE-based controller synthesis method. Numerical simulations are shown in Section~\ref{section:5}, and we conclude the paper in Section~\ref{section:6}.

\emph{Notations}: The fields of real numbers and $m\times n$ real matrices are denoted by $\mathbb{R}$ and $\mathbb{R}^{m\times n}$, respectively. The closed right-half complex plane is denoted by $\bar{\mathbb{C}}^+$ . A matrix $ M\in\mathbb{R}^{m\times n}$ is represented by its entries $m_{ij},i ={1,2,\dots,m},j = {1,2,\dots,n},$ for convenience, \emph{i.e.}, $M=[m_{ij}]$, and its transpose is denoted by $M^T$. The spectrum of a square matrix $M\in\mathbb{R}^{n\times n}$ is denoted by $\sigma(M)$. A matrix is called Hurwitz (or stable) if and only if all of its eigenvalues have negative real parts. An $n\times n$ diagonal matrix, whose diagonal entries are $m_1,m_2,\dots,m_n$ and start from the upper left, is denoted by $\text{diag}\{m_1,m_2,\dots,m_n\}$. For matrices $A=[a_{ij}]\in\mathbb{R}^{m\times n}$ and $B\in\mathbb{R}^{p\times q}$, the Kronecker product of $A$ and $B$ is denoted by $A\otimes B$. For any positive integer $N$, let $\mathcal{N}=\{1,2,\dots,N\}$, and the identity matrix of dimension $N$ is denoted by $I_N$.

\section{Problem Statement: Cooperative Control of Heterogeneous Connected Vehicles} \label{section:2}

\begin{figure}[!t]
    \centering
    \includegraphics[width=0.9\columnwidth]{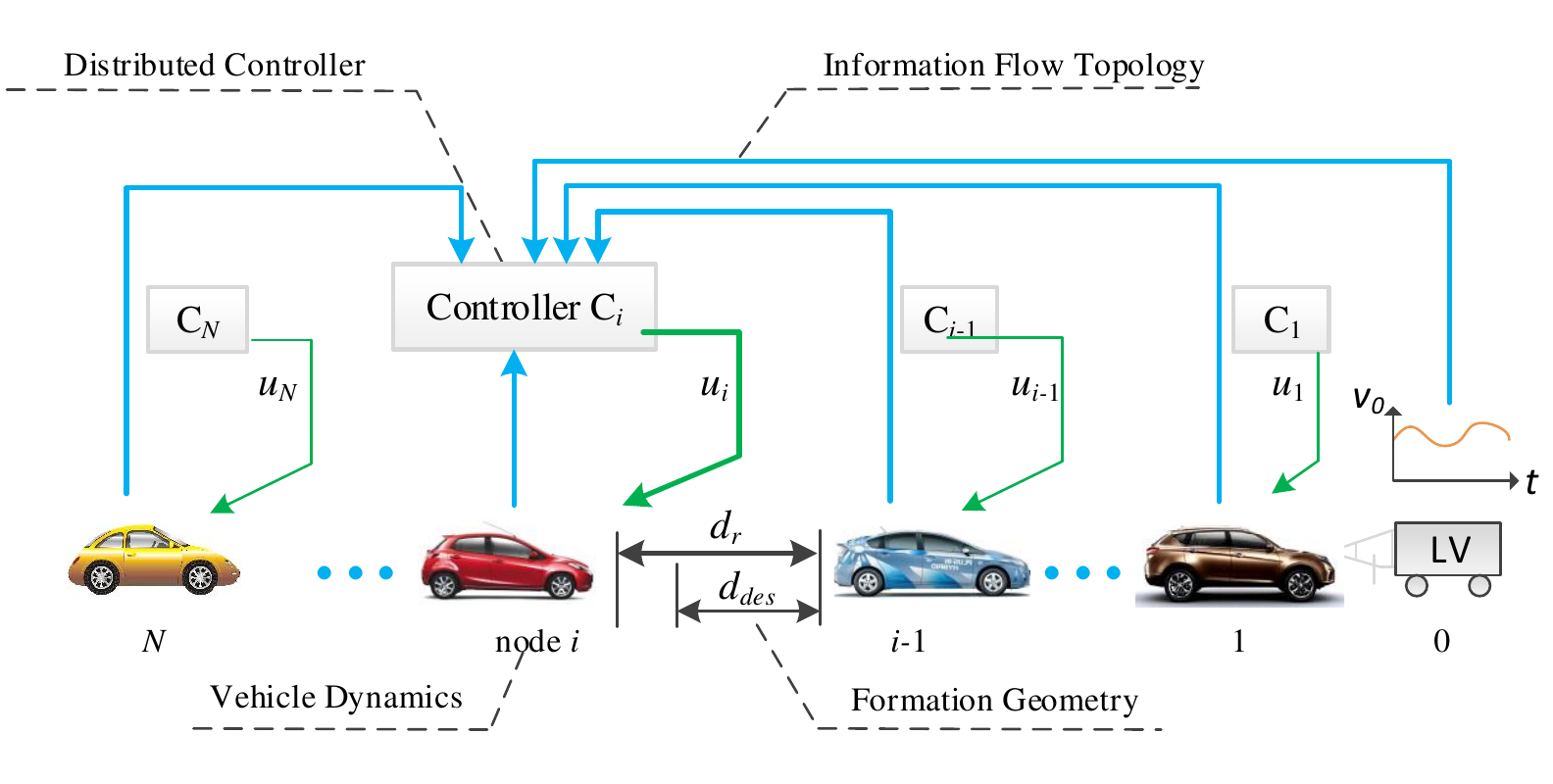}
    \caption{Four major components of a heterogeneous vehicle platoon~\cite{zheng2015dynamic}.} 
    \label{fig:1}
\end{figure}

We consider the cooperative control of a platoon of $N+1$ heterogeneous vehicles (nodes) running on a straight flat road, consisting of a leading vehicle indexed by $0$ and $N$ following vehicles indexed from $1$ to $N$ (see Fig.~\ref{fig:1}). The control objective is to make the following vehicles move at the same velocity as the leading vehicle while maintaining a fixed formation geometry. As shown in Fig.~\ref{fig:1}, from a control perspective, the platoon can be viewed as a combination of four main components: 1) vehicle dynamics; 2) distributed controller; 3) information flow topology; 4) formation geometry, which is known as the four-component framework~\cite{zheng2015dynamic}. A categorization of platoon control can be found in~\cite{zheng2015dynamic} and~\cite{li2015overview} based on the features of each component.

In this section, we briefly introduce the modeling of the four components and present the problem statement of platooning of connected vehicles.

\subsection{Modeling of Heterogeneous Platoons}
\subsubsection{Vehicular longitudinal dynamics}

Vehicle dynamics describe the behavior of each node, which are inherently nonlinear, consisting of the engine, brake system, and rolling resistance, \emph{etc.} However, a detailed nonlinear model may not be helpful for platoon level analysis, since it in general cannot lead to analytical results, especially when we consider the effect of different communication topologies. In the literature, to strike a balance between accuracy and conciseness, a typical choice is to derive a linear model by either using a hierarchical control framework~\cite{shladover1991automated}, or employing a feedback linearization technique~\cite{zheng2016stability}. The linear model is then served as a basis for theoretical analysis for the cooperative control of a platoon.

Here, we use the following linearized third-order model to describe the longitudinal behavior of each vehicle in a platoon
\begin{eqnarray}\label{eq:1}
  \begin{cases}
  \dot{p}_i(t)=v_i(t), \\
  \dot{v}_i(t)=a_i(t), \\
  \tau_i\dot{a}_i(t)+a_i(t)=u_i(t),
  \end{cases} \forall i\in \{0\}\cup \mathcal{N},
\end{eqnarray}
where $p_i(t)$, $v_i(t)$, and $a_i(t)$ denote the position, velocity, and acceleration of the $i$-th node, respectively; $u_i(t)$ represents the desired acceleration of the $i$-th node, and $\tau_i$ characterizes the inertial time lag of the powertrain system. This model is widely used as a basis for platoon level analysis; see,~\emph{e.g.},~\cite{li2017dynamical,zheng2015dynamic,ploeg2014controller,ploeg2014lp}. Here, we note that the third equation in~\eqref{eq:1} is a first-order inertial function that approximates the acceleration response of vehicle longitudinal dynamics. The parameter $\tau_i$ can be different, \emph{e.g.}, $\tau_i$ is small for passenger cars while it is big for commercial cars. Concisely, Eq.~\eqref{eq:1} can be rewritten as 
\begin{eqnarray}\label{eq:2}
  \begin{cases}
  \dot{x}_i(t) = A_ix_i(t)+ B_iu_i(t), \\
  y_i(t) = Cx_i(t),
  \end{cases}
\end{eqnarray}
where $x_i =\begin{bmatrix}
        p_i,
        v_i,
        a_i
        \end{bmatrix}^T$,
\begin{equation*}
    \begin{aligned}
    A_i=\begin{bmatrix}
          0 & 1 & 0 \\
          0 & 0 & 1 \\
          0 & 0 & \displaystyle -\frac{1}{\tau_i} \\
        \end{bmatrix},
        B_i=\begin{bmatrix}
          0 \\
          0 \\
          \displaystyle \frac{1}{\tau_i}  \\
        \end{bmatrix},C=\begin{bmatrix}
          c_p & 0 & 0 \\
          0 & c_v & 0 \\
          0 & 0 & c_a \\
        \end{bmatrix},
        \end{aligned}
\end{equation*}
and $c_\sharp\in\{0,1\}$, $\sharp=p$, $v$, or $a$, denotes whether the corresponding state of the vehicle can be measured as an output.

\subsubsection{Model for information flow topologies}

Information flow topologies describe how the vehicles in a platoon exchange information with each other, which exerts a great influence on the collective behavior of a platoon. We employ directed graphs to model the information flow topology in a platoon; please refer to~\cite{godsil2013algebraic} for more details on graph theory. In this study, we assume the communication is perfect and ignore the effects such as data quantization and time delay.

The information flow between followers is modeled by a directed graph $\mathcal{G} (\mathcal{V},\mathcal{E})$ with a set of nodes $\mathcal{V} =\{v_1,v_2,\ldots,v_N\}$ representing each following vehicle and a set of edges $\mathcal{E}\subseteq \mathcal{V} \times \mathcal{V} $ representing the information exchange. The adjacency matrix $\mathcal{A}$ associated with $\mathcal{G}$ is defined as $\mathcal{A}=[a_{ij}]\in\mathbb{R}^{N\times N}$ with each entry
\begin{equation*}
  a_{ij}=\begin{cases}
            1& \text{if }\{v_j,v_i\}\in\mathcal{E},\\
            0& \text{otherwise},
        \end{cases}
\end{equation*}
where $\{v_j,v_i\}\in\mathcal{E}$ means vehicle $i$ can obtain the information of vehicle $j$. It is assumed that there is no self-loop, \emph{i.e.}, $a_{ii}=0$. The degree matrix $\mathcal{D}$ is defined as $\mathcal{D}=\text{diag}\{d_{11},d_{22},\dots,d_{NN}\}\in\mathbb{R}^{N\times N}$, where
\begin{equation*}
    d_{ii}=\sum_{k=1}^N a_{ik}.
\end{equation*}
The Laplacian matrix $\mathcal{L}$ associated with $\mathcal{G}$ is defined as $\mathcal{L}=[l_{ij}]\in \mathbb{R}^{N\times N}$ with each entry:
\begin{equation}\label{eq:5}
    l_{ij}=\begin{cases}
            -a_{ij} & i\neq j\\
            \sum_{k=1}^N a_{ik} & i=j
            \end{cases}.
\end{equation}
Then, we have
\begin{equation*}
\mathcal{L}=\mathcal{D}-\mathcal{A}.
\vspace{-2mm}
\end{equation*}

To model the connections between the leader and followers, we define a pinning matrix:
\begin{equation}\label{eq:6}
  \mathcal{P}=\text{diag}\{p_{11},p_{22},\dots,p_{NN}\},
\end{equation}
where $p_{ii}=1$ means node $i$ can obtain the information from the leader; otherwise $p_{ii}=0$. Then, the connections in a platoon are described by the matrices $\mathcal{L}$ and $\mathcal{P}$, which are naturally suitable for different communication topologies.

\subsubsection{Formation geometry}

Formation geometry describes the desired inter-vehicle spacing between two adjacent nodes, which is the main objective in the cooperative control of a platoon. Mathematically, we require
\begin{equation}\label{eq:7}
  \begin{cases}
  \displaystyle
  \lim_{t\to+\infty}\left\|v_i(t)-v_0(t)\right\|=0,\\
  \displaystyle
  \lim_{t\to+\infty}\left\|p_i(t)-p_{i-1}(t)-d_{\text{des},i-1,i}\right\|=0,
  \end{cases} \forall i\in \mathcal{N},
\end{equation}
where $d_{\text{des},i-1,i}$ is the desired gap between vehicle $i-1$ and vehicle $i$. This value can be either velocity-dependent (referred to as the constant time headway policy), or velocity- independent (called the constant spacing policy). In this paper, we use the constant spacing policy, \emph{i.e.}, $d_{\text{des},i-1,i}=-d_0$, as used in~\cite{barooah2009mistuning,zheng2015stabilityMargin,herman2015nonzero}.

\subsubsection{Design of distributed controllers}

The distributed controller defines a feedback law using local information that is available for each node, such that the collective behavior of a platoon reaches the global coordination~\eqref{eq:7}. In this paper, we consider a linear feedback of the form
\begin{equation}\label{eq:8}
  u_i=-k_i^T\sum_{j=1}^N a_{ij}(y_i-y_j-\hat{d}_{i,j})-p_{ii}k_i^T(y_i-y_0-\hat{d}_{i,0}),
\end{equation}
where $k_i=\begin{bmatrix} k_{ip}, k_{iv}, k_{ia}\end{bmatrix}^T$ is the local feedback gain and $\hat{d}_{i,j}=\begin{bmatrix} (i-j)d_0, 0, 0 \end{bmatrix}^T$ is the desired distance vector. In Eq.~\eqref{eq:8}, only local information is used for feedback. We note that the controller~\eqref{eq:8} and its variants have been widely used in~\cite{zheng2016stability,barooah2009mistuning,herman2015nonzero,zheng2015stabilityMargin}.

\subsection{Problem Statement of Heterogeneous Platoon Design}
There are two heterogeneous sources in a practical platoon: 1) heterogeneous dynamics $(A_i,B_i)$; 2) heterogeneous feedback gains $k_i$. In previous studies~\cite{gao2016robust,li2012distributed}, the dynamical heterogeneity is considered as a type of uncertainty, \emph{i.e.}, $A_i=A+\Delta A_i,B_i=B+\Delta B_i$, where the nominal system $(A,B)$ is homogeneous and the uncertainty is bounded. Note that the authors of~\cite{gao2016robust,li2012distributed} still assume homogenous feedback gains, \emph{i.e.}, $k_i=k$, even in the case of heterogeneous dynamics.

In this paper, we directly consider the heterogeneity in dynamics $(A_i,B_i)$ and allow for heterogeneous feedback gains $k_i$ as well. This treatment not only offers more freedom in controller design, but also gains a deeper understanding of the influence of heterogeneity on the collective behavior of a platoon. Nevertheless, this treatment brings new challenges for the analysis and design since we cannot use the decomposition results in the homogeneous case and new tools are needed. Precisely, we seek to address the following issues:
\begin{enumerate}
  \item \textit{Feasibility issue}: We address whether it is feasible to achieve the goal \eqref{eq:7} using a controller of the form \eqref{eq:8} for heterogeneous platoons.
  \item \textit{Stability region}: We derive an analytical region of $k_i$ where the closed-loop platoon is asymptotically stable.
  \item \textit{Controller synthesis}: We introduce a method to calculate a particular $k_i$ that stabilizes the closed-loop platoon.
\end{enumerate}

These three issues are related to each other. Unlike the homogeneous case, the feasibility of the controller \eqref{eq:8} becomes nontrivial due to the heterogeneity. It requires careful discussions on the existence of the controller \eqref{eq:8} that can achieve the goal \eqref{eq:7}. In this paper, we apply the internal model principle of multi-agent systems~\cite{wieland2011internal} to heterogeneous platoon design, which explicitly highlights tracking ability of each following vehicle. Then, an internal stability theorem is derived using the Routh-Hurwitz stability criterion by exploiting a lower-triangular structure. We further propose an ARE based design method to calculate the feedback gain $ k_i$ for each following vehicle in a heterogeneous platoon.

\section{Tracking Ability of Heterogeneous Platoons}\label{section:3}
In this section, we formally address the feasibility of the controller~\eqref{eq:8} by using the internal model principle of heterogeneous multi-agent systems.

\begin{myLem}[Internal model principle~\cite{wieland2011internal,seyboth2015robust}] \label{lemma:1}
     Consider a heterogeneous linear network of $N$ agents \eqref{eq:2} with static diffusive couplings \eqref{eq:8}. If $\lim_{t\to+\infty}\|y_i(t)-y_j(t)-\hat{d}_{i,j}\|=0, \forall i,j\in\{0\}\cup\mathcal{N}$, then there exists an integer $m>0$ and three matrices $\Pi_i\in \mathbb{R}^{3\times m}$ with full column rank, $S\in \mathbb{R}^{m\times m}$ and $R\in \mathbb{R}^{3\times m}$, where $\sigma(S)\in\bar{\mathbb{C}}^+$ and $(S,R)$ is observable, such that
     \begin{subequations}
        \begin{align}
            A_i\Pi_i&=\Pi_iS, \label{eq:9}\\
            C\Pi_i&=R.  \label{eq:10}
        \end{align}
     \end{subequations}
\end{myLem}

This result is applicable for general heterogeneous multi-agent systems, known as the \emph{internal model principle}, which is originally proved in~\cite{wieland2011internal}. A special form of this theorem appeared in~\cite{seyboth2015robust}. In principle, the conditions in Lemma \ref{lemma:1} indicate that all followers are able to track the leader defined by the dynamical matrix $S$ and the output matrix $R$. Also, the dynamics of the followers must embed an internal model of the leader; the interested reader is referred to~\cite{wieland2011internal,seyboth2015robust} for more details.

In terms of heterogeneous platoons, we have the following result.

\begin{myTheo}\label{theorem:1}
    Consider the cooperative control of a platoon of heterogeneous connected vehicles with dynamics defined in~\eqref{eq:2} and controller given by~\eqref{eq:8}. If the cooperation objective~\eqref{eq:7} is satisfied, then the leader must move at a constant speed, \emph{i.e.}, $\dot{p}_0=v_0, \dot{v}_0=0$, and the spacing information must be measured for each vehicle, \emph{i.e.}, $c_p=1$.
\end{myTheo}

\begin{IEEEproof}
    According to Lemma \ref{lemma:1}, the conditions in Eq. \eqref{eq:9} and Eq. \eqref{eq:10} must be satisfied. Since $\Pi_i$ has full column rank, we know that each eigenvalue of $S$ is also an eigenvalue of $A_i$, \emph{i.e.},
    \begin{equation}\label{eq:11}
      \sigma (S)\subseteq\sigma (A_i),\forall i\in \mathcal{N}.
    \end{equation}
    Then, the eigenvalues of $S$ are a subset of the largest common set $\displaystyle \cap_{i=1}^N \sigma(A_i)$. According to the dynamics~\eqref{eq:2}, we have
    $$\bigcap_{i=1}^N \sigma(A_i)=\{0\}.$$
    Therefore, the matrix $S$ can only have zero eigenvalues. In this case, possible choices for $S$ include
    $$
        S_1 = 0, S_2 = \begin{bmatrix}
                0 & 0 \\
                0 & 0
                \end{bmatrix}, S_3 = \begin{bmatrix}
                0 & 1 \\
                0 & 0
                \end{bmatrix}.
    $$

    The first two choices lead to a trivial situation $\dot{p}_0 = v_0 =0$, where the leader's velocity is zero. This is a special case where the leader move at a constant speed. Here, we consider a broader case: for heterogeneous platoons, the dynamical matrix $S$ of the leader is
    \begin{equation}\label{eq:12}
      S=\begin{bmatrix}
                0 & 1 \\
                0 & 0
                \end{bmatrix},
    \end{equation}
    which indicates $\dot p_0=v_0,\dot v_0=0$. In this case, the leader can have a non-zero constant velocity. Also, we know that the following matrices $\Pi_i$ and $R$ satisfy the conditions in Eq. \eqref{eq:9} and Eq.~\eqref{eq:10}
    \begin{equation}\label{eq:13}
      \Pi_i=\begin{bmatrix}
                    1 & 0 \\
                    0 & 1 \\
                    0 & 0
            \end{bmatrix},
      R=\begin{bmatrix}
                    c_p & 0 \\
                    0 & c_v \\
                    0 & 0
                  \end{bmatrix}.
    \end{equation}
    The observability matrix of $(S,R)$ is
    \begin{equation}\label{eq:14}
      Q_o=\begin{bmatrix}
        c_p & 0 \\
        0 & c_v \\
        0 & 0 \\
        0 & c_p \\
        0 & 0 \\
        0 & 0\end{bmatrix}.
    \end{equation}

    Then, the observability of $(S,R)$ requires $\text{rank}(Q_o)=2$, implying $c_p\neq 0$. In our case, it means the spacing information is available to each vehicle, \emph{i.e.}, $c_p=1$.
\end{IEEEproof}

In a platoon of connected vehicles, the leader's state actually defines the equilibrium point of each follower, \emph{i.e.}, each follower tries to reach consensus on the equilibrium state defined by the leader (implicitly or explicitly). In Theorem~\ref{theorem:1}, we formally show that to reach the steady consensus state, the leader should run at a constant speed, implying that the followers can only track at most the leader's spacing and velocity and that there should be no acceleration in the leader. Note that this is the objective defined in~\eqref{eq:7}.

In fact, many previous studies directly assume that the leader's velocity is constant; see, \emph{e.g.},~\cite{dunbar2012distributed,Zheng2016distributed,zheng2016stability,zheng2015stabilityMargin,barooah2009mistuning,herman2015nonzero}. The transient from one constant speed to another is usually modeled as a certain disturbance of the leader, where the response of each follower is typically studied using the notion of string stability~\cite{shaw2007string,naus2010string,ploeg2014lp}. In addition, the conditions in Theorem~\ref{theorem:1} also explicitly highlight the necessity of spacing information for controller feedback.

\begin{myRem}\label{remark:1}
    For homogenous platoons, a feasible solution to Eq. \eqref{eq:9} and Eq. \eqref{eq:10} is trivial, \emph{i.e.}, $S=A,R=C,\Pi_i=I_3$, because homogeneous followers share the same internal model with the leader. We note that the conditions in Theorem \ref{theorem:1} only highlight the necessary requirements on the dynamics. For controllability, the communication graph  should contain at least one spanning tree rooting at the leader~\cite{ren2005consensus}, \emph{i.e.}, there should exist a directed path from the leader to every follower (a \emph{pinning condition} is required). In other words, the leader's information should be available to every follower explicitly or implicitly. As we shall see below, this requirement is confirmed in Theorem \ref{theorem:2}.
\end{myRem}

\section{Stability Region And Controller Synthesis}\label{section:4}

The last section gives some necessary conditions for the platoon design. With these conditions in mind, this section discusses the stability region of heterogeneous platoons with directed acyclic interactions. Also, a synthesis method is proposed based on the solution to an ARE.

\subsection{Closed-loop Platoon Dynamics}
Here, we first formulate the closed-loop dynamics of heterogeneous platoons. The desired trajectory of the follower $i$ is shown as
\begin{equation*}
  p_i=p_0-i\times d_0.
\end{equation*}
As shown in Theorem \ref{theorem:1}, we assume the leader moves at a constant speed, \emph{i.e.},
$\dot p_0=v_0, \ddot p_0=0.$
We then define the following tracking error for each follower:
\begin{equation}\label{eq:16}
    \begin{cases}
        \hat{p}_i=p_i-p_0+i\times d_0, \\
        \hat{v}_i=\dot{p}_i-\dot{p}_0=v_i-v_0, \\
        \hat{a}_i=\ddot{p}_i-\ddot{p}_0=a_i.
    \end{cases}
\end{equation}
Further, the lumped tracking error of follower $i$ is
\begin{equation}\label{eq:17}
  \epsilon_i=\sum_{j=1}^N a_{ij}(\hat{y}_i-\hat{y}_j)+p_{ii}\hat{y}_i,
\end{equation}
where $\hat{y}_i=C\hat{x}_i$ and $\hat{x}_i=[\hat{p}_i,\hat{v}_i,\hat{a}_i]^T$ denotes the output and state of tracking errors, respectively. Then, the control law~\eqref{eq:8} can be rewritten into a compact form
\begin{equation}\label{eq:18}
u_i=-k_i^T\epsilon_i.
\end{equation}

The closed-loop dynamics of tracking errors are written as
\begin{equation}\label{eq:19}
    \begin{cases}
        \dot{\hat{p}}_i = \hat{v}_i, \\
        \dot{\hat{v}}_i = \hat{a}_i, \\
        \displaystyle \dot{\hat{a}}_i = -\frac{1}{\tau_i}\hat{a}_i-\frac{1}{\tau_i}k_i^T\epsilon_i.
    \end{cases}
\end{equation}
Also, we know that
\begin{equation*}
  \frac{1}{\tau_i}k_i^T\epsilon_i=\sum_{j=1}^Na_{ij}\left(\frac{1}{\tau_i}k_i^TC(\hat{x}_i-\hat{x}_j)\right)+p_{ii}\frac{1}{\tau_i}k_i^TC\hat{x}_i.
\end{equation*}
For simplicity, we define
\begin{equation} \label{eq:20}
    \begin{aligned}
         \frac{1}{\tau_i}k_i^TC &=   \begin{bmatrix} \displaystyle \frac{1}{\tau_i}k_{ip}c_{p} &  \displaystyle \frac{1}{\tau_i}k_{iv}c_{v} & \displaystyle \frac{1}{\tau_i}k_{ia}c_{a} \end{bmatrix}\\
    &\equiv \begin{bmatrix}t_{ip} & t_{iv} & t_{ia}\end{bmatrix},
    \end{aligned}
\end{equation}
where a set of new variables are introduced,
$$t_{ip} =  \frac{1}{\tau_i}k_{ip}c_{p}, t_{iv} =  \frac{1}{\tau_i}k_{iv}c_{v}, t_{ia} =  \frac{1}{\tau_i}k_{ia}c_{a}.$$

Then, the closed-loop platoon dynamics can be compactly rewritten as
\begin{equation}\label{eq:21}
  \begin{bmatrix}
          \dot{\hat p} \\
          \dot{\hat v} \\
          \dot{\hat a}
        \end{bmatrix}=
  \begin{bmatrix}
            0 & I_N & 0 \\
            0 & 0 & I_N \\
            -T_pG & -T_vG & -\Delta - T_aG \\
        \end{bmatrix}
  \begin{bmatrix}
          \hat p \\
          \hat v \\
          \hat a
        \end{bmatrix},
\end{equation}
where $G=\mathcal{L+P}$, and
\begin{equation*}
  \hat p=\begin{bmatrix}
          \hat p_1 \\
          \hat p_2 \\
          \vdots \\
          \hat p_N
          \end{bmatrix},
  \hat v=\begin{bmatrix}
          \hat v_1 \\
          \hat v_2 \\
          \vdots \\
          \hat v_N\end{bmatrix},
  \hat a=\begin{bmatrix}
          \hat a_1 \\
          \hat a_2 \\
          \vdots \\
          \hat a_N\end{bmatrix}
\end{equation*}
denote the lumped states of tracking errors, and
\begin{equation}\label{eq:22}
    \Delta=\begin{bmatrix}
                \displaystyle \frac{1}{\tau_1} &  &  \\
                 & \ddots &  \\
                 &  & \displaystyle \frac{1}{\tau_N} \\
              \end{bmatrix}
\end{equation}
collects the effects of heterogeneous inertial time lags of followers, and
\begin{equation}\label{eq:23}
    T_\sharp=\begin{bmatrix}
                t_{1\sharp} &  &  \\
                 & \ddots & \\
                 &  & t_{N\sharp} \\
              \end{bmatrix},\sharp \in \{p,v,a\},
\end{equation}
assembles the effects of heterogeneous feedback gains.

\begin{myRem}\label{remark:2}
    In case of homogeneous platoons (\emph{i.e.}, $A_i=A,B_i=B,k_i=k,\forall i\in\mathcal{N})$, $x_i$ is usually used as the state variable. As shown in~\cite{zheng2016stability}, the closed-loop dynamics of homogeneous platoons can be written as
    \begin{equation}\label{eq:24}
        \dot{x}=(I_N\otimes A-G\otimes Bk^T)x,	
    \end{equation}
    where $x=[x_1^T,x_2^T,\dots,x_N^T]^T$. However, when it comes to heterogeneous platoons, the closed-loop dynamics will be more complicated than~\eqref{eq:24} when choosing $x$ as the state variable. In contrast, we collect $\hat p$, $\hat v$, and $\hat a$ as the state variables to construct the closed-loop platoon dynamics, leading to the concise formulation \eqref{eq:21}. In~\eqref{eq:21}, both the heterogeneous inertial time lags $\tau_i$ and controller gains $[k_{ip},k_{iv},k_{ia}]$ are collected together into diagonal matrices (see Eq.~\eqref{eq:22} and Eq.~\eqref{eq:23}), making it easier to analyze their effects on the closed-loop system. In the homogeneous case, the matrices $\Delta$ and $T_\sharp$ become homogeneous as well in a form of $\eta I_N$, and~\eqref{eq:21} can be transformed into~\eqref{eq:24} using a certain state transformation.
\end{myRem}

\begin{myRem}\label{remark:3}
    From~\eqref{eq:21}, it is easy to see that each of the four components in Fig. \ref{fig:1} exerts a certain influence on the platoon dynamics: the matrix $\Delta$ and the structure of the state matrix represent the longitudinal dynamics of vehicles; the matrix $G$ represents the influence of information flow topology; the vectors $\hat p$, $\hat v$, and $\hat a$ contain the effect of formation geometry; and the matrix $T_\sharp$ shows the effect of the distributed controller. This is consistent with the homogeneous cases; see~\eqref{eq:24}.
\end{myRem}

\begin{figure}[t]
  \centering
    \includegraphics[width=0.95\columnwidth]{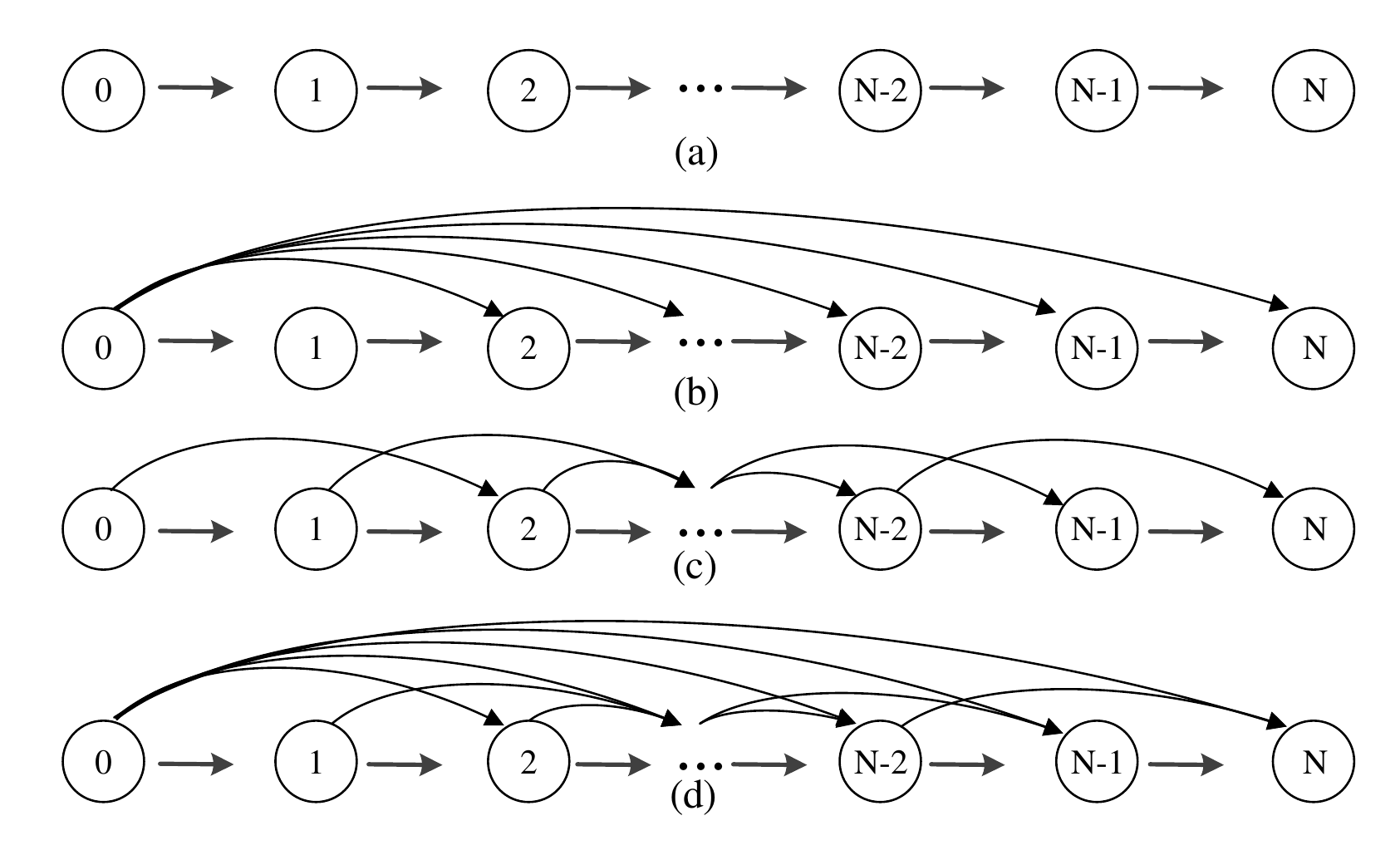}
  \caption{Examples of unidirectional topologies: (a) PF, (b) PLF, (c) two-predecessor following (TPF), (d) two-predecessor leader following (TPLF).}
  \label{fig:3}
\end{figure}

\subsection{Directed Acyclic Graph}\label{section:4b}

We now introduce the definition of directed acyclic graphs.
\begin{myDef}\label{def:1}
    A directed acyclic graph (DAG) is a finite directed graph with no directed cycles.
\end{myDef}

Equivalently, a DAG is a directed graph that has a topological ordering, \emph{i.e.}, a sequence of the vertices, such that every edge is directed unidirectionally from preceding nodes to downstream ones in the sequence. In fact, the DAG is a natural extension of unidirectional topologies (see \emph{e.g.},~\cite{Zheng2016distributed}), including the common PF~\cite{ploeg2014controller} and PLF~\cite{naus2010string} topologies as special cases (see Fig.~\ref{fig:3}).

\begin{figure}[t]
    \centering
    \setlength{\abovecaptionskip}{1em}
    \footnotesize
    \begin{tikzpicture}
	  \matrix (m) [matrix of nodes,
	  		       row sep = 0.8em,	
	  		       column sep = 1.6em,	
  			       nodes={circle, draw=black}] at (-2,0)
         {  1 & 2 & 3 & 4 \\};
		\draw[-latex] (m-1-1) -- (m-1-2);
		\draw[-latex] (m-1-2) -- (m-1-3);
        \draw[-latex] (m-1-4) -- (m-1-3);
        \draw[-latex] (m-1-1)  to [out=45,in=150] (m-1-4);
        \draw[-latex] (m-1-1)  to [out=-45,in=-150] (m-1-3);
		\node at (-2,-0.8) {(a)};

		\matrix (m2) [matrix of nodes,
	  		       row sep = 0.8em,	
	  		       column sep = 1.6 em,	
  			       nodes={circle, draw=black}] at (2,0)
         {  1 & 2 & 4 & 3 \\};
		\draw[-latex] (m2-1-1) -- (m2-1-2);
		\draw[-latex] (m2-1-2) to [out=45,in=150] (m2-1-4);
        \draw[-latex] (m2-1-3) -- (m2-1-4);
        \draw[-latex] (m2-1-1)  to [out=45,in=135] (m2-1-4);
        \draw[-latex] (m2-1-1)  to [out=-45,in=-150] (m2-1-3);
        \node at (2,-0.8) {(b)};

        \matrix (m3) [matrix of nodes,
	  		       row sep = 0.8em,	
	  		       column sep = 1.6em,	
  			       nodes={circle, draw=black}] at (-2,-1.8)
         {  1 & 4 & 2 & 3 \\};
		\draw[-latex] (m3-1-1) -- (m3-1-2);
		\draw[-latex] (m3-1-2) to [out=45,in=150] (m3-1-4);
        \draw[-latex] (m3-1-3) -- (m3-1-4);
        \draw[-latex] (m3-1-1)  to [out=45,in=135] (m3-1-4);
        \draw[-latex] (m3-1-1)  to [out=-45,in=-150] (m3-1-3);
		\node at (-2,-2.6) {(c)};

		\matrix (m4) [matrix of nodes,
	  		       row sep = 0.8em,	
	  		       column sep = 1.6 em,	
  			       nodes={circle, draw=black}] at (2,-1.8)
         {  1 & 2 & 3 & 4 \\};
		\draw[-latex] (m4-1-1) -- (m4-1-2);
		\draw[-latex] (m4-1-2) -- (m4-1-3);
        \draw[-latex] (m4-1-4) -- (m4-1-3);
        \draw[-latex] (m4-1-1)  to [out=45,in=150] (m4-1-4);
        \draw[-latex] (m4-1-3)  to [out=-135,in=-30] (m4-1-1);
        \node at (2,-2.6) {(d)};
	\end{tikzpicture}
\vspace{-1mm}
    \caption{Examples of DAGs and topological ordering: (a)-(c) are DAGs while (d) is not.}
    \label{fig:2}
\end{figure}
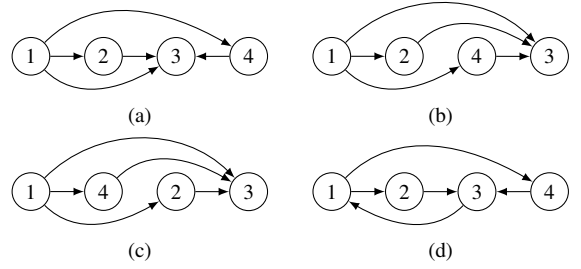

For example, as shown in Fig.~\ref{fig:2}, graphs (a)-(c) are DAGs, because they are all directed graphs with no directed cycles. Also, Fig.~\ref{fig:2}(b) and Fig.~\ref{fig:2}(c) are two types of topological ordering for Fig.~\ref{fig:2}(a), where the order of vertices is shifted such that only information flow from preceding nodes (the left side) to downstream nodes (the right side) is allowed. On the contrary, the graph in Fig.~\ref{fig:2}(d), which can be obtained by reversing the information flow direction between nodes $1$ and $3$ in Fig.~\ref{fig:2}(a), is not a DAG, since nodes $\{1,2,3\}$ and nodes $\{1,4,3\}$ form two directed cycles.

Based on Definition \ref{def:1}, the following lemma gives a formulation of the topological ordering.

\begin{myLem}\label{lemma:2}
    For a DAG with a set of nodes indexed as $\mathcal{N}=\{1,2,\dots,N\}$, there exists at least one permutation, denoted by the ordered set $\{s_1,s_2,\dots,s_N\}$, such that
    \begin{equation}\label{eq:25}
        \hat a_{ij}=0,\forall i<j,i,j\in\mathcal{N},
    \end{equation}
    where $\hat a_{ij}$ is the entry of the adjacency matrix associated with the permutated graph. Equivalently, there exists an invertible permutation matrix $Q\in \mathbb{R}^{N\times N}=[e_{s_1},e_{s_2},\dots,e_{s_N}]$, where $e_i\in \mathbb{R}^N$ is the standard unit vector in the $i$-th direction,  such that
    \begin{subequations}
      \begin{align}
      \begin{bmatrix}
            s_1 \\
            s_2 \\
            \vdots \\
            s_N
            \end{bmatrix}=Q
      \begin{bmatrix}
            1 \\
            2 \\
            \vdots \\
            N
            \end{bmatrix}, \label{eq:26}\\
      \hat{\mathcal{A}}=Q^{-1}\mathcal{A}Q, \label{eq:27}
      \end{align}
    \end{subequations}
    where $\mathcal{A}$ and $\hat{\mathcal{A}}$ are the adjacency matrices of the original and permuted graphs, and $\hat{\mathcal{A}}$ is a lower-triangular matrix.
\end{myLem}

\begin{floatEq}
    \begin{equation}\label{eq:29}
    \begin{aligned}
    |\lambda I_{3N}-\hat{A}|&=\left|
            \begin{array}{ccc}
              \lambda I_N & -I_N & 0 \\
              0 & \lambda I_N & -I_N \\
              T_pG & T_vG & \lambda I_N+\Delta+T_aG
            \end{array}
         \right|\cdot
         \left|
            \begin{array}{ccc}
              I_N & \frac{1}{\lambda}I_N & \frac{1}{\lambda^2}I_N \\
              0 & I_N & \frac{1}{\lambda}I_N \\
              0 & 0 & I_N
            \end{array}
         \right| \\
         &=\left|
            \begin{array}{ccc}
              \lambda I_N & 0 & 0 \\
              0 & \lambda I_N & 0 \\
              T_pG & \displaystyle \frac{1}{\lambda}T_pG+T_vG & \displaystyle \frac{1}{\lambda^2}T_pG+\frac{1}{\lambda}T_vG+\lambda I_N+\Delta+T_aG
            \end{array}
         \right| \\
         &=|\lambda^3 I_N+\lambda^2(\Delta+T_aG)+\lambda T_vG+T_pG|.
         \end{aligned}
         \end{equation}
\end{floatEq}

\begin{floatEq}
\begin{equation}\label{eq:30}
    \begin{aligned}
    |\lambda I_{3N}-\hat{A}|
    &=|Q^{-1}||\lambda^3 I_N+\lambda^2(\Delta+T_aG)+\lambda T_vG+T_pG||Q| \\
    &=|\lambda^3 I_N+\lambda^2(\hat\Delta+\hat{T}_a\hat{G})+\lambda \hat{T}_v\hat{G}+\hat{T}_p\hat{G}| \\
    &=\prod_{i=1}^N\left(\lambda^3+\lambda^2\left(\frac{1}{\tau_i}+t_{ia}(d_{ii}+p_{ii})\right)+\lambda t_{iv}(d_{ii}+p_{ii})+t_{ip}(d_{ii}+p_{ii})\right).
    \end{aligned}
\end{equation}
\end{floatEq}

The proof is straightforward, and we omit it for brevity. Take graphs (a) and (b) in Fig. \ref{fig:2} for example. The ordered vertex sets of graphs (a) and (b) are $\{1,2,3,4\}$ and $\{1,2,4,3\}$, respectively, and the adjacency matrices are
\begin{equation*}
  \mathcal{A}_a=\begin{bmatrix}
                    0 & 0 & 0 & 0 \\
                    1 & 0 & 0 & 0 \\
                    1 & 1 & 0 & 1 \\
                    1 & 0 & 0 & 0
                   \end{bmatrix},
  \mathcal{A}_b=\begin{bmatrix}
                    0 & 0 & 0 & 0 \\
                    1 & 0 & 0 & 0 \\
                    1 & 0 & 0 & 0 \\
                    1 & 1 & 1 & 0
                    \end{bmatrix},
\end{equation*}
respectively. Then we know graph (b) is a permutation of graph (a) with the invertible permutation matrix
\begin{equation*}
  Q=\begin{bmatrix}
                    1 & 0 & 0 & 0 \\
                    0 & 1 & 0 & 0 \\
                    0 & 0 & 0 & 1 \\
                    0 & 0 & 1 & 0
      \end{bmatrix},
\end{equation*}
which satisfies
\begin{equation*}
  \begin{bmatrix}
        s_1 \\
        s_2 \\
        s_3 \\
        s_4
        \end{bmatrix}=
  \begin{bmatrix}
        1 \\
        2 \\
        4 \\
        3
        \end{bmatrix}=Q
  \begin{bmatrix}
        1 \\
        2 \\
        3 \\
        4
        \end{bmatrix},
\end{equation*}
  \begin{equation*}
\mathcal{A}_{b}=Q^{-1}\mathcal{A}_{a}Q,
\end{equation*}
and $\mathcal{A}_b$ is a lower-triangular matrix.

The key fact of Lemma~\ref{lemma:2} is that the adjacency matrix of a DAG can be transformed into a lower-triangular matrix, which facilitates the analysis and synthesis of heterogeneous platoons subsequently. Furthermore, for the permutation matrix $Q$ in Lemma~\ref{lemma:2}, we have the following result.

\begin{myLem}\label{lemma:3}
    Consider a matrix $Q\in\mathbb{R}^{N\times N}=[e_{s_1},e_{s_2},\dots,e_{s_N}]$, where $e_i\in\mathbb{R}^N$ is the
standard unit vector in the $i$-th direction, and $\{s_1,s_2,\dots,s_N\}$ is a permutation of $\{1,2,\dots,N\}$. Then, for any diagonal matrix $F\in\mathbb{R}^{N\times N}$, the matrix $\hat F=Q^{-1}FQ$ is still diagonal with the same diagonal entries as $F$, but the order of the entries is permuted.
\end{myLem}

The proof is straightforward since $\hat F=Q^{-1}FQ$ represents certain column and row operations on the diagonal matrix $F$. Thus, $\hat F=Q^{-1}FQ$ remains diagonal, and only the order of its diagonal entries is permuted according to the permutation $\{s_1,s_2,\dots,s_N\}$.

According to Eq.~\eqref{eq:27} in Lemma~\ref{lemma:2} and Lemma~\ref{lemma:3}, one direct result is that for DAGs, the associated matrix $G=\mathcal{L+P}=\mathcal{D-A+P}$ can be transformed into a lower-triangular matrix $\hat G=Q^{-1}GQ$ after permutation, since
$$
    Q^{-1}GQ=Q^{-1}\mathcal{(D+P)}Q-\hat{\mathcal{A}},
$$
and $\mathcal{D}$ and $\mathcal{P}$ are diagonal matrices which remain diagonal after transformation.

\subsection{Stability Region Analysis}
We are ready to present the second theorem of this paper.

\begin{myTheo}\label{theorem:2}
    Consider the cooperative control of a platoon of heterogeneous connected vehicles with the closed-loop dynamics given by~\eqref{eq:21}. If the information flow topology is a DAG, then the platoon is asymptotically stable if and only if the following statements hold.
    \vspace{-1.1mm}
\begin{enumerate}
  \item The spacing and velocity information are measurable, \emph{i.e.}, $c_p=c_v=1$;
  \item Every follower can obtain the information of at least one other node, \emph{i.e.}, $d_{ii}+p_{ii}>0$;
  \item The local feedback gains satisfy
    \begin{equation}\label{eq:28}
      \begin{cases}
        k_{ip}>0, \\
        \displaystyle
        k_{iv}>\frac{\tau_ik_{ip}}{1+k_{ia}c_a(d_{ii}+p_{ii})}, \\
        \displaystyle
        k_{ia}c_a>-\frac{1}{d_{ii}+p_{ii}},
      \end{cases} \forall i\in\mathcal{N}.
    \end{equation}
\end{enumerate}
\end{myTheo}
\begin{IEEEproof}
For simplicity, in Eq.~\eqref{eq:21}, we denote
\begin{equation*}
  \hat A=\begin{bmatrix}
              0 & I_N & 0 \\
              0 & 0 & I_N \\
              -T_pG & -T_vG & -\Delta-T_aG
           \end{bmatrix}.
\end{equation*}
We first analyze the characteristic equation of the closed-loop platoon system \eqref{eq:21}:
\begin{equation*}
  |\lambda I_{3N}-\hat{A}|=\left|
            \begin{array}{ccc}
              \lambda I_N & -I_N & 0 \\
              0 & \lambda I_N & -I_N \\
              T_pG & T_vG & \lambda I_N+\Delta+T_aG
            \end{array}
         \right|.
\end{equation*}
If the closed-loop system is asymptotically stable, then there are no zero roots, \emph{i.e.}, $\lambda_i\neq 0,\forall i\in\mathcal{N}$. Consequently, we have~\eqref{eq:29}.

If the information flow topology is a DAG, then according to Lemma \ref{lemma:3}, we know that there exists a permutation matrix $Q\in \mathbb{R}^{N\times N}$ such that $\hat G=Q^{-1}GQ$ is a lower-triangular matrix. Then, we denote
\begin{equation*}
    \begin{cases}
        \hat\Delta =Q^{-1}\Delta Q, \\
         \hat T_\sharp=Q^{-1}T_\sharp Q,\sharp\in\{p,v,a\},
    \end{cases}
\end{equation*}
which are diagonal matrices and have the same diagonal entries as $\Delta$ and $T_\sharp$ according to Lemma \ref{lemma:3}. These facts lead to the equation shown in~\eqref{eq:30}.

Therefore, the stability of the closed-loop platoon system \eqref{eq:21} is equivalent to the stability of the following $N$ characteristic equations:
\begin{multline}\label{eq:31}
  \lambda^3+\lambda^2\left(\frac{1}{\tau_i}+t_{ia}(d_{ii}+p_{ii})\right) \\
  +\lambda t_{iv}(d_{ii}+p_{ii})+t_{ip}(d_{ii}+p_{ii})=0,i\in\mathcal{N}.
\end{multline}
The stability of~\eqref{eq:31} can then be checked by using the Routh-Hurwitz stability criterion, as shown in~\eqref{eq:32}.
\begin{equation}\label{eq:32}
  \begin{array}{ccc}
    \lambda^3 & 1 & t_{iv}(d_{ii}+p_{ii}) \\
    \lambda^2 & \displaystyle \frac{1}{\tau_i}+t_{ia}(d_{ii}+p_{ii}) & t_{ip}(d_{ii}+p_{ii}) \\
    \lambda^1 & \displaystyle t_{iv}(d_{ii}+p_{ii})-\frac{t_{ip}(d_{ii}+p_{ii})}{\frac{1}{\tau_i}+t_{ia}(d_{ii}+p_{ii})} &  \\
    \lambda^0 & t_{ip}(d_{ii}+p_{ii}) &
  \end{array}
\end{equation}
Thus, we have
\begin{equation}\label{eq:33}
  \begin{cases}
  \displaystyle
    \frac{1}{\tau_i}+t_{ia}(d_{ii}+p_{ii})>0, \\
    \displaystyle
    t_{iv}(d_{ii}+p_{ii})-\frac{t_{ip}(d_{ii}+p_{ii})}{\frac{1}{\tau_i}+t_{ia}(d_{ii}+p_{ii})}>0, \\
    t_{ip}(d_{ii}+p_{ii})>0.
  \end{cases}
\end{equation}

Then, it is easy to know that $t_{ip}\neq0,t_{iv}\neq0$, which indicates $c_p\neq0,c_v\neq0$. In our case, it requires $c_p=c_v=1$, \emph{i.e.}, the spacing and velocity information are measurable. Also, from the third inequality of~\eqref{eq:33}, we know
\begin{equation*}
  d_{ii}+p_{ii}>0.
\end{equation*}
After some simple linear algebra, we arrive at the requirements on the local feedback gains, shown in Eq.~\eqref{eq:28}. This completes the proof.

\end{IEEEproof}

There are a number of points that are worth highlighting for Theorem~\ref{theorem:2}. The first condition in Theorem~\ref{theorem:2} is consistent with Theorem~\ref{theorem:1}, and it states that both spacing and velocity information are necessary to stabilize a heterogeneous platoon system. Also, this result agrees with the earliest platooning practices~\cite{shladover1991automated}, where only spacing and velocity information are available since the sensing systems are often radar-based and lack the acceleration information of other vehicles. In this case, the condition~\eqref{eq:28} in Theorem~\ref{theorem:2} reduces to
\begin{equation}\label{eq:34}
  k_{ip}>0,k_{iv}>\tau_ik_{ip},
\end{equation}
which guarantees the internal stability of a heterogeneous platoon without relying on the acceleration information. The second condition in Theorem \ref{theorem:2} means there should exist at least one spanning tree rooting at the leader, which agrees with~\cite{ren2005consensus} (see Remark~\ref{remark:1}). It is easy to check that all the unidirectional topologies shown in Fig. \ref{fig:3} satisfy this property. In addition, the condition~\eqref{eq:28} generalizes the stability condition in~\cite{zheng2016stability} to heterogeneous platoons with directed acyclic interactions.

\begin{myRem}\label{remark:4}
    It is shown that the acceleration information via V2V communication helps improve string stability~\cite{oncu2014cooperative,ploeg2014controller}. For condition~\eqref{eq:28}, it is easy to see that the feedback gain $k_{ia}$ actually enlarges the stability region of $k_{iv}$, whose lower bound is reduced from $\tau_ik_{ip}$ to $$\frac{\tau_ik_{ip}}{1+k_{ia}c_a(d_{ii}+p_{ii})}.$$
    The additional freedom brought by acceleration feedback could then be used to improve other performance indexes, \emph{e.g.}, string stability. From this perspective, our result is consistent with the statement in~\cite{oncu2014cooperative,ploeg2014controller} as well.
\end{myRem}

\subsection{ARE-based Design of Feedback Gains}

Theorem \ref{theorem:2} gives analytical results on the stability region, within which all the feedback gains in Eq.~\eqref{eq:28} can guarantee asymptotical stability of the system~\eqref{eq:21}. However, Theorem~\ref{theorem:2} does not indicate how to choose a proper control gain for a specific platoon system.

In this section, based on Theorem~\ref{theorem:2}, we present a feedback gain design method according to the solution to an ARE. This method has a relatively clear physical interpretation, which is easy to use in practice.

\begin{myTheo}\label{theorem:3}
    Consider the cooperative control of a platoon of heterogeneous connected vehicles with closed-loop dynamics given by~\eqref{eq:21} and all states measurable, \emph{i.e.}, $c_p=c_v=c_a=1$. If the information flow topology is a DAG with a spanning tree, \emph{i.e.}, $d_{ii}+p_{ii}>0$, then the control gain given in~\eqref{eq:35}, where $P_i\succeq 0$ is the root of the ARE~\eqref{eq:36}, guarantees the stability of system~\eqref{eq:21} if $\displaystyle \alpha_i\geq\frac{1}{2(d_{ii}+p_{ii})}, \forall i\in\mathcal{N}$.
    \begin{subequations}
        \begin{align}
            k_i^T &=\alpha_iB_i^TP_i, \label{eq:35}\\
            P_iA_i+A_i^TP_i&-P_iB_iB_i^TP_i+\varepsilon_iI_3=0,\varepsilon_i>0. \label{eq:36}
        \end{align}
    \end{subequations}
\end{myTheo}
\begin{IEEEproof}
    According to~\eqref{eq:30}, we first analyze the stability of the matrix
    \begin{equation}\label{eq:42}
        \tilde A_i=A_i-(d_{ii}+p_{ii})B_ik_i^T.
    \end{equation}
    Consider a Lyapunov equation using the positive definite solution $P_i$ to the ARE~\eqref{eq:36}:
    \begin{equation}\label{eq:37}
      P_i\tilde A_i+\tilde A_i^TP_i=-\tilde{Q}_i.
    \end{equation}
    When $\displaystyle \alpha_i\geq\frac{1}{2(d_{ii}+p_{ii})}, \forall i\in\mathcal{N}$, by substituting Eq. \eqref{eq:35} into Eq.~\eqref{eq:36}, we have
    \begin{equation*}
        \begin{aligned}
      \tilde{Q}_i&=-\left(1-2\alpha_i(d_{ii}+p_{ii})\right)P_iB_iB_i^TP_i+\varepsilon_iI_3 \\
      &\succeq \varepsilon_iI_3 \succ 0.
      \end{aligned}
    \end{equation*}
    Then, the matrix $\tilde A_i, \forall i\in\mathcal{N}$ is Hurwitz. This means all the eigenvalues of its characteristic equation, shown in Eq.~\eqref{eq:38}, have negative real parts:
    \begin{multline}\label{eq:38}
      |\lambda I_{3N}-\tilde{A}_i|=\lambda^3+\lambda^2\left(\frac{1}{\tau_i}+t_{ia}(d_{ii}+p_{ii})\right) \\
        +\lambda t_{iv}(d_{ii}+p_{ii})+t_{ip}(d_{ii}+p_{ii}) = 0 .
    \end{multline}

    Therefore, according to~\eqref{eq:31} in Theorem \ref{theorem:2}, the closed-loop system \eqref{eq:21} is stable.

    \end{IEEEproof}

Compared to the condition~\eqref{eq:28} in Theorem \ref{theorem:2}, where three gain parameters should lie within an explicit and analytic region, Theorem~\ref{theorem:3} gives an implicit and univariate way to choose feedback gains. 
The ARE~\eqref{eq:36} corresponds to the infinite time LQR problem for system \eqref{eq:2} with a performance index $$\displaystyle J_i=\frac{1}{2}\int_0^\infty(\varepsilon_ix_i^Tx_i+u_i^2)dt.$$ Since $(A_i,B_i)$ is controllable, the ARE~\eqref{eq:36} always has a unique positive definite solution $P_i$ for any $\varepsilon_i>0$.

\begin{myRem}\label{remark:5}
    The equation \eqref{eq:37} actually defines a Lyapunov function $V_i=x_i^TP_ix_i$ with derivative $\dot V_i=-x_i^T\tilde Q_ix_i$, for the subsystem $\dot x_i=\tilde A_ix$, where $\tilde A_i$ is given in~\eqref{eq:42}. Then, a larger $\varepsilon_i$ implies a faster convergence rate of this subsystem. In addition, according to Eq.~\eqref{eq:30} and Eq.~\eqref{eq:38}, the relative stability of this subsystem is equivalent to that of the original platoon system. This fact brings much convenience to adjust the convergence rate of the platoon system by properly choosing the parameter $\varepsilon_i$. Note that the control gains given by the ARE-based method is only a subset of all the stable control gains given in Theorem \ref{theorem:2}; see Eq.~\eqref{eq:28}.
\end{myRem}

\begin{myRem}\label{remark:6}
    Note that the stability of the closed-loop system is equivalent to the stability of the subsystem $\dot x_i=\tilde A_ix$, where $\tilde A_i$ is given in ~\eqref{eq:42}. We have used this fact in the proofs of Theorems~\ref{theorem:1} and~\ref{theorem:2}. In homogenous cases with general information flow topologies, it is shown in~\cite{fax2004information} that the equivalent system becomes $\tilde A_i=A_i-\lambda_iB_ik_i^T$, where $\lambda_i$ is also the eigenvalue of the matrix $G=\mathcal{L+P}$.
\end{myRem}

\section{Numerical Results}\label{section:5}

This section presents numerical simulations to validate the effectiveness of our findings. In particular, we consider a heterogeneous platoon with eight vehicles (one leading vehicle and seven following vehicles) under multiple types of directed acyclic interaction topologies, including PF, PLF, TPF, and TPLF topologies (see Fig.~\ref{fig:3}). Also, simulations with a realistic nonlinear vehicle model are carried out to show effectiveness of our results in real traffic environments.

\subsection{Validations based on the linear model}
In the simulations, the desired spacing was set to $d_0=20~m$. The initial states of the leading vehicle were $p_0(0)=0~m, v_0(0)=20~m/s$ with the velocity profile given by
\begin{equation}\label{eq:39}
  v_0(t)=\begin{cases}
            10 & 0 s\leq t<3 s \\
            10+t & 3 s\leq t<15 s \\
            22 & t\geq15 s
        \end{cases},(m/s).
\end{equation}
The initial states of all following vehicles were set as $p_i (0)=-i\times d_0,v_i (0)=v_0 (0),a_i (0)=0,\forall i\in \mathcal{N}$. We assume that the states of each vehicle are all measurable, \emph{i.e.}, $c_p=c_v=c_a=1$.

\begin{table}[t]
    \centering
    \renewcommand\arraystretch{1.2}
    \caption{The time lag $\tau_i$ and feedback gains of each follower.}
    \begin{tabular}{ m{1.1cm}<{\centering}| m{0.8cm}<{\centering} m{0.8cm}<{\centering} m{0.8cm}<{\centering} m{0.8cm}<{\centering} m{0.8cm}<{\centering}}
        \hline \toprule[1pt]
        Vehicle Index & $\tau_i[s]$ & $k_{ip}$ & $k_{iv}$ & $\hat{k}_{iv}$ & $k_{ia}$ \\
        \hline
        1 & 0.40 & 3.00 & 3.40 & 0.06 & 2.00 \\
        2 & 0.55 & 1.30 & 3.55 & 0.09 & 2.62 \\
        3 & 0.32 & 2.31 & 3.32 & 0.10 & 2.87 \\
        4 & 0.44 & 1.65 & 3.44 & 0.08 & 2.97 \\
        5 & 0.38 & 3.83 & 3.38 & 0.07 & 3.07 \\
        6 & 0.51 & 2.42 & 3.51 & 0.05 & 3.70 \\
        7 & 0.29 & 2.91 & 3.29 & 0.04 & 2.79 \\
        \bottomrule[1pt]
    \end{tabular}
    \label{tab:1}
\end{table}

\begin{figure}[t]
  \centering
    \setlength{\belowcaptionskip}{0em}
  \subfigure[ ]{\includegraphics[width=0.45\columnwidth]{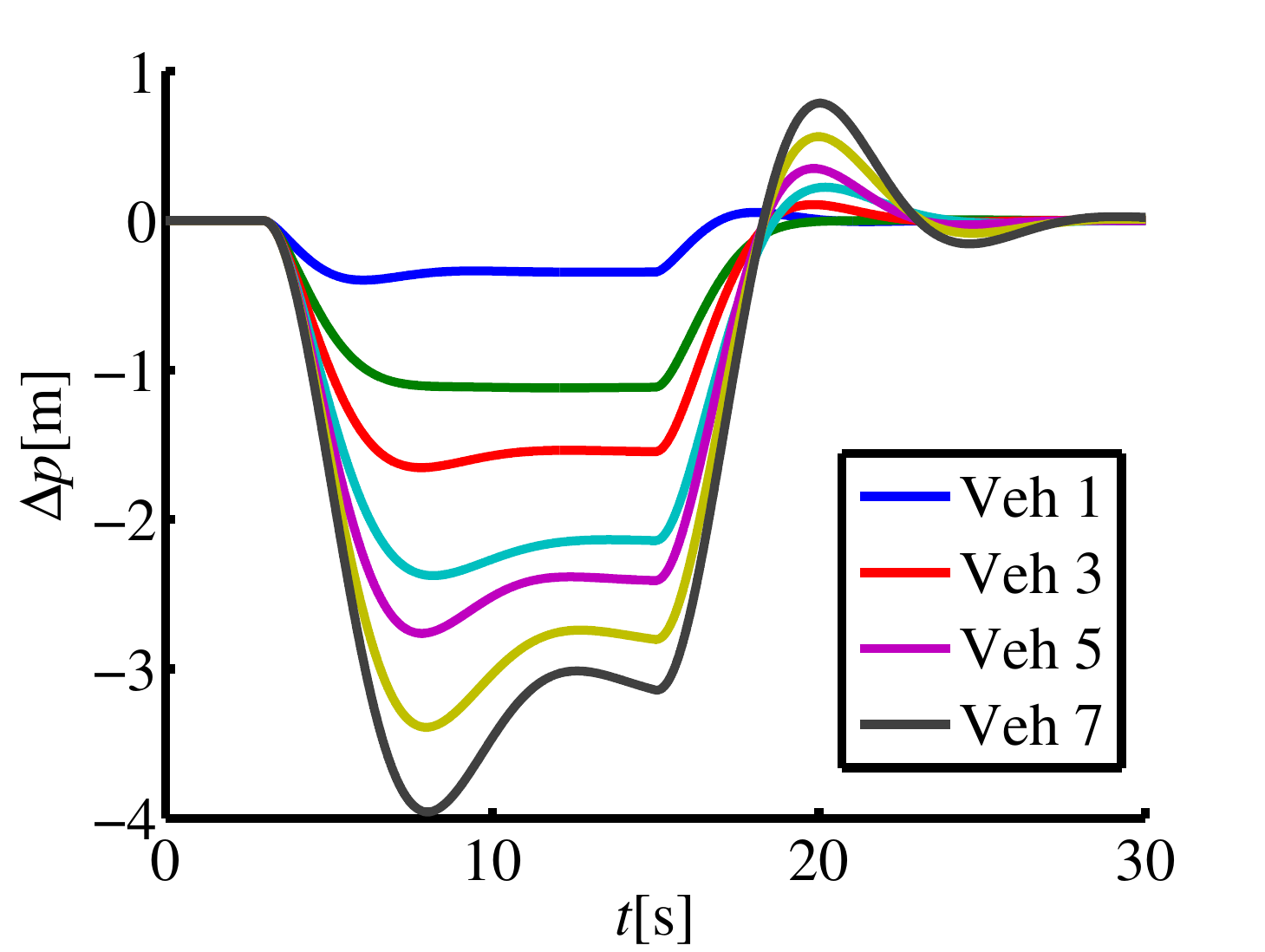}} \hspace{0.1 cm}
  \subfigure[ ]{\includegraphics[width=0.45\columnwidth]{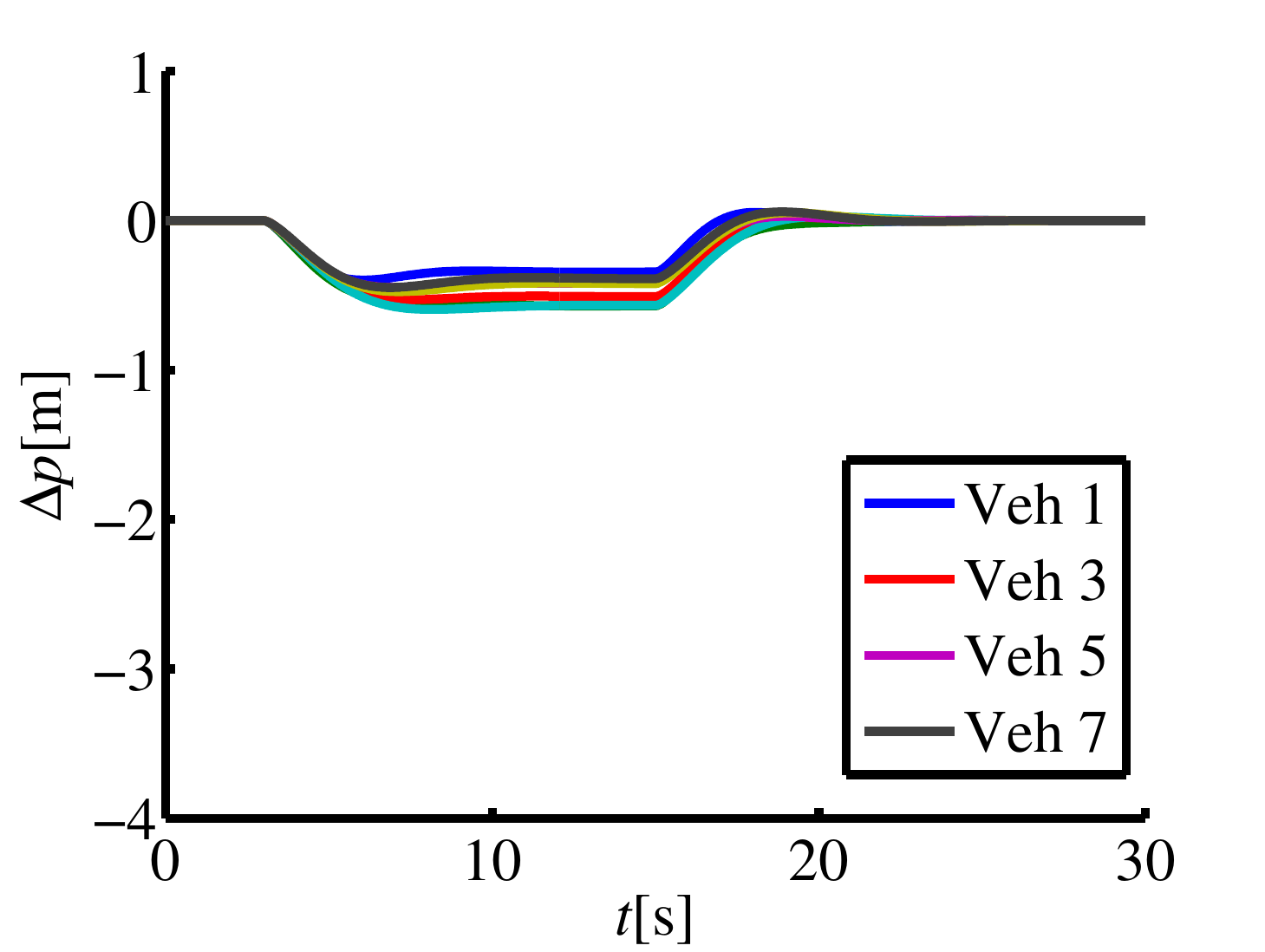}} 
  \\ \vspace{-2mm}
  \subfigure[ ]{\includegraphics[width=0.45\columnwidth]{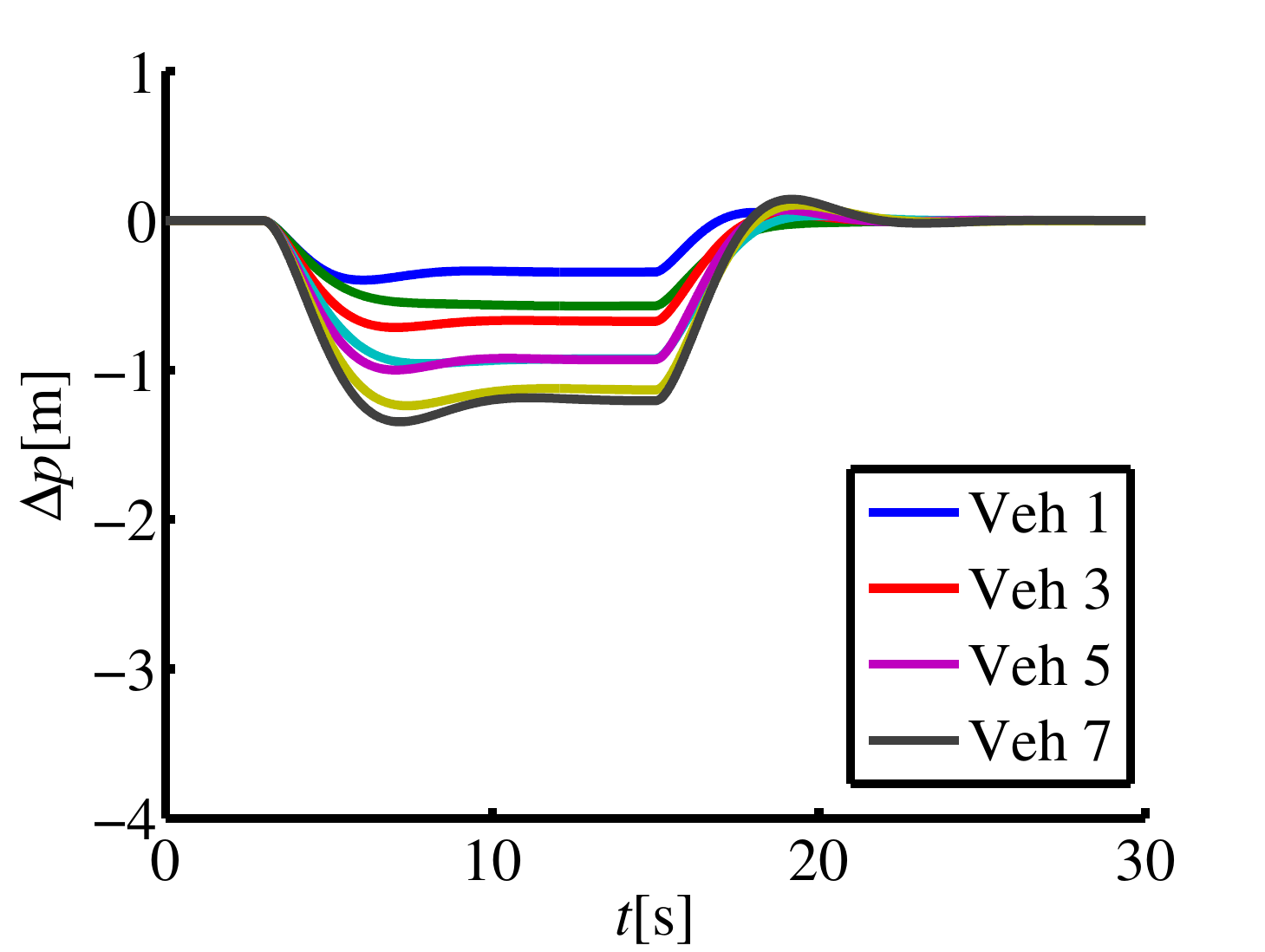}} \hspace{0.1 cm}
  \subfigure[ ]{\includegraphics[width=0.45\columnwidth]{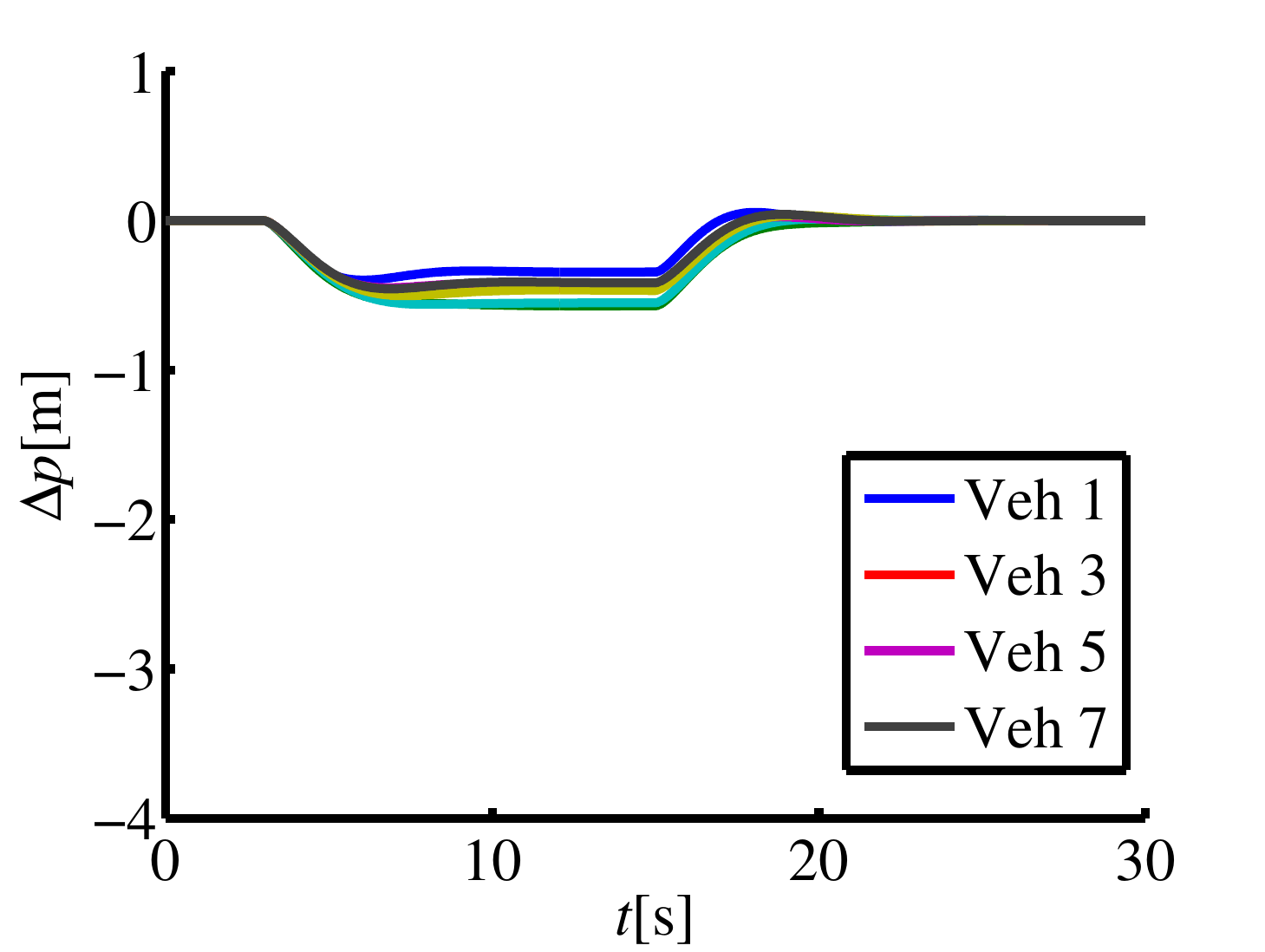}}
  \caption{Spacing error profiles using the feedback gains $(k_{ip},k_{iv},k_{ia})$ in Table~\ref{tab:1} (Theorem \ref{theorem:2} is satisfied) and the velocity profile \eqref{eq:39} for heterogeneous platoons with different topologies: (a) PF, (b) PLF, (c) TPF, (d) TPLF.}
  \label{fig:4}
\end{figure}

\begin{figure}[t]
  \centering
  \subfigure[ ]{\includegraphics[width=0.45\columnwidth]{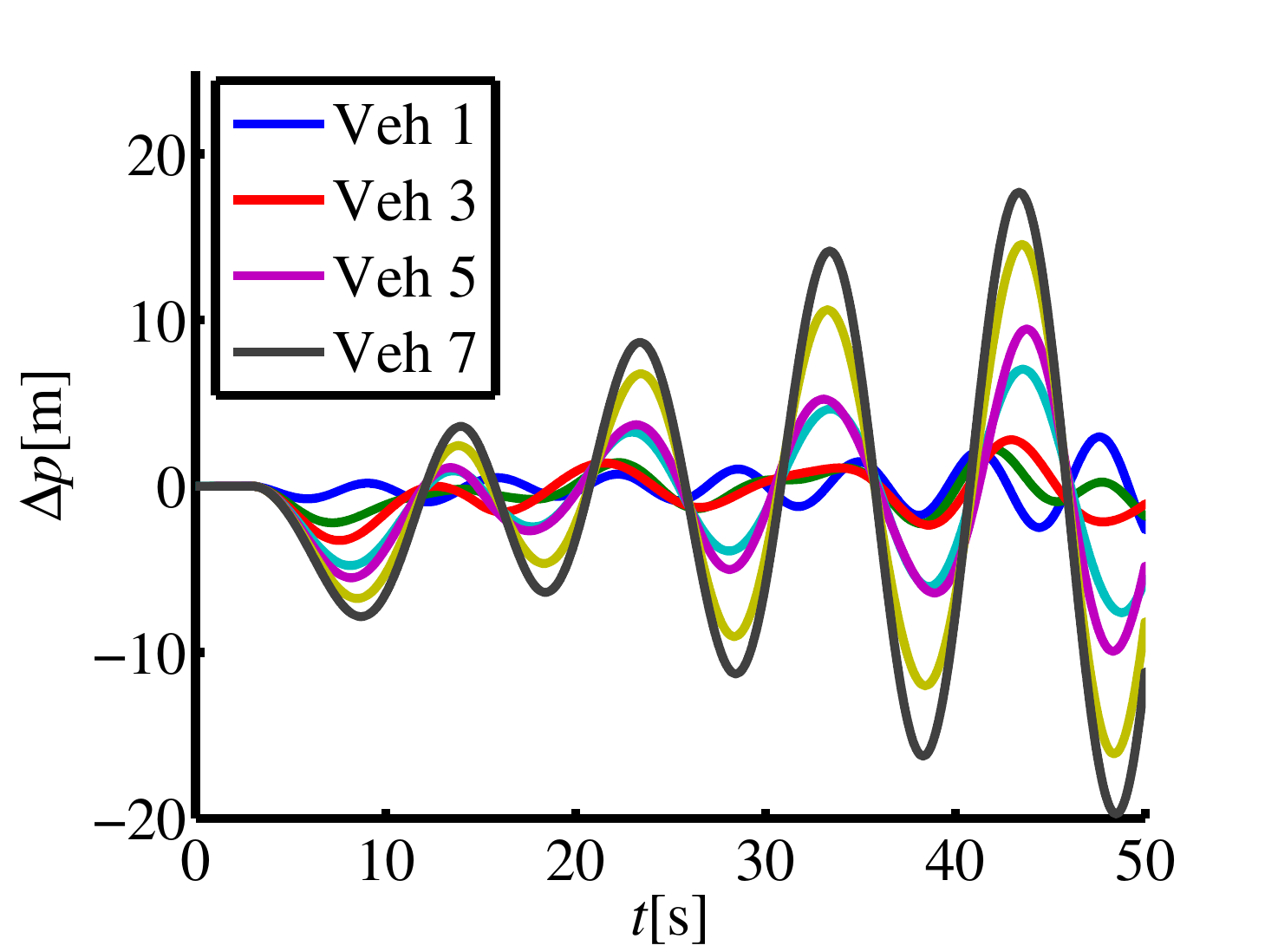}}
  \subfigure[ ]{\includegraphics[width=0.45\columnwidth]{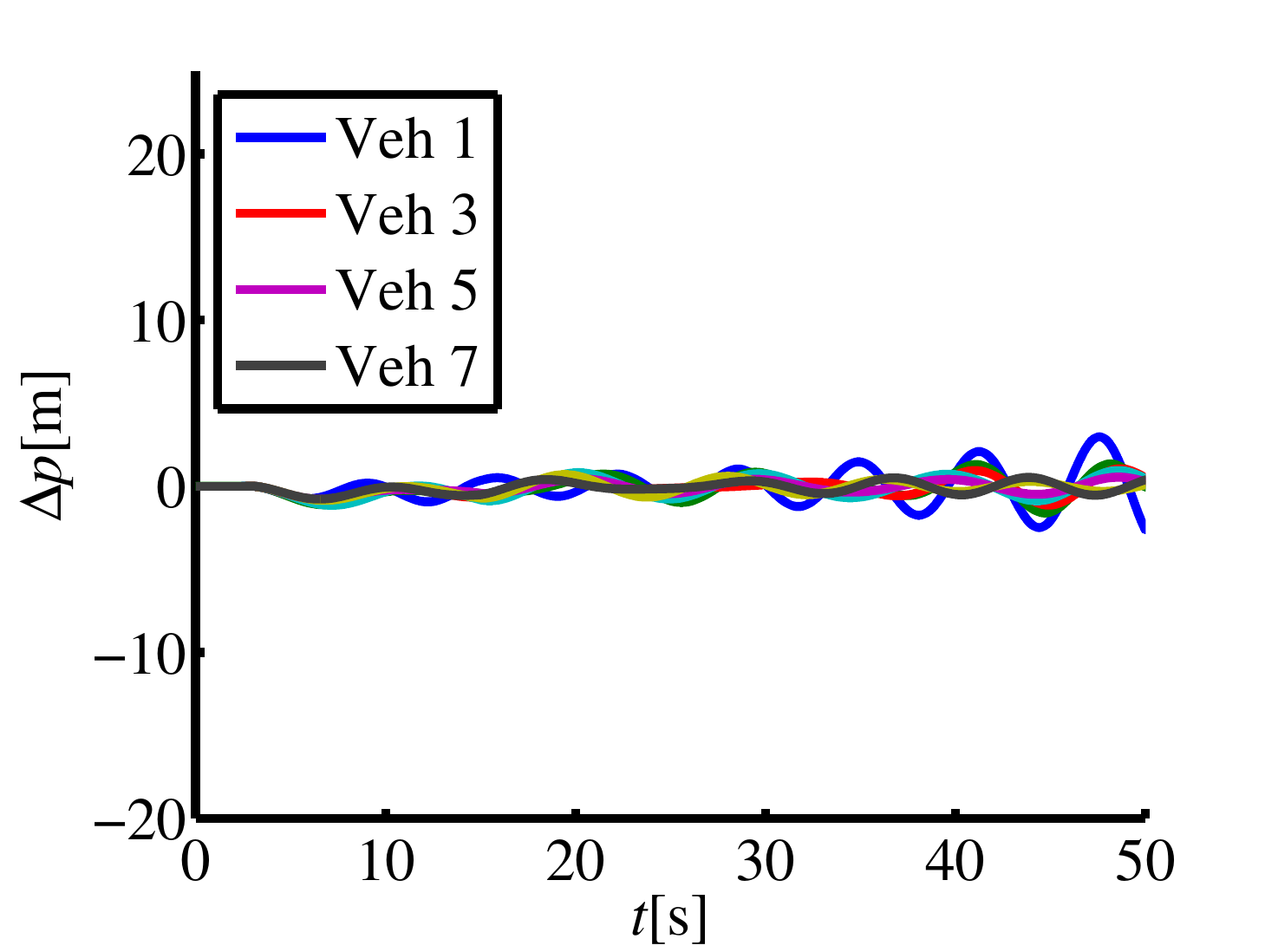}}
  \subfigure[ ]{\includegraphics[width=0.45\columnwidth]{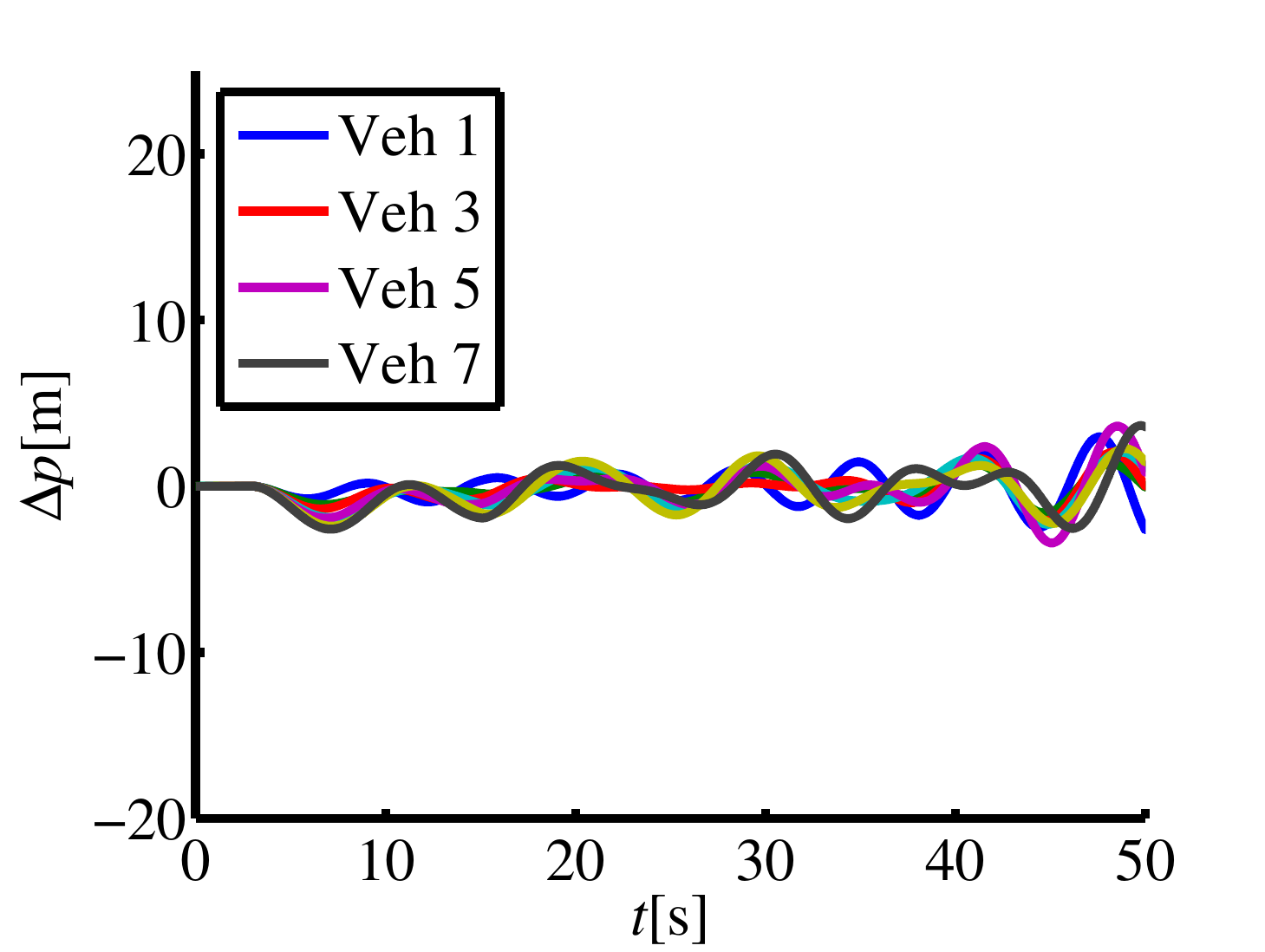}}
  \subfigure[ ]{\includegraphics[width=0.45\columnwidth]{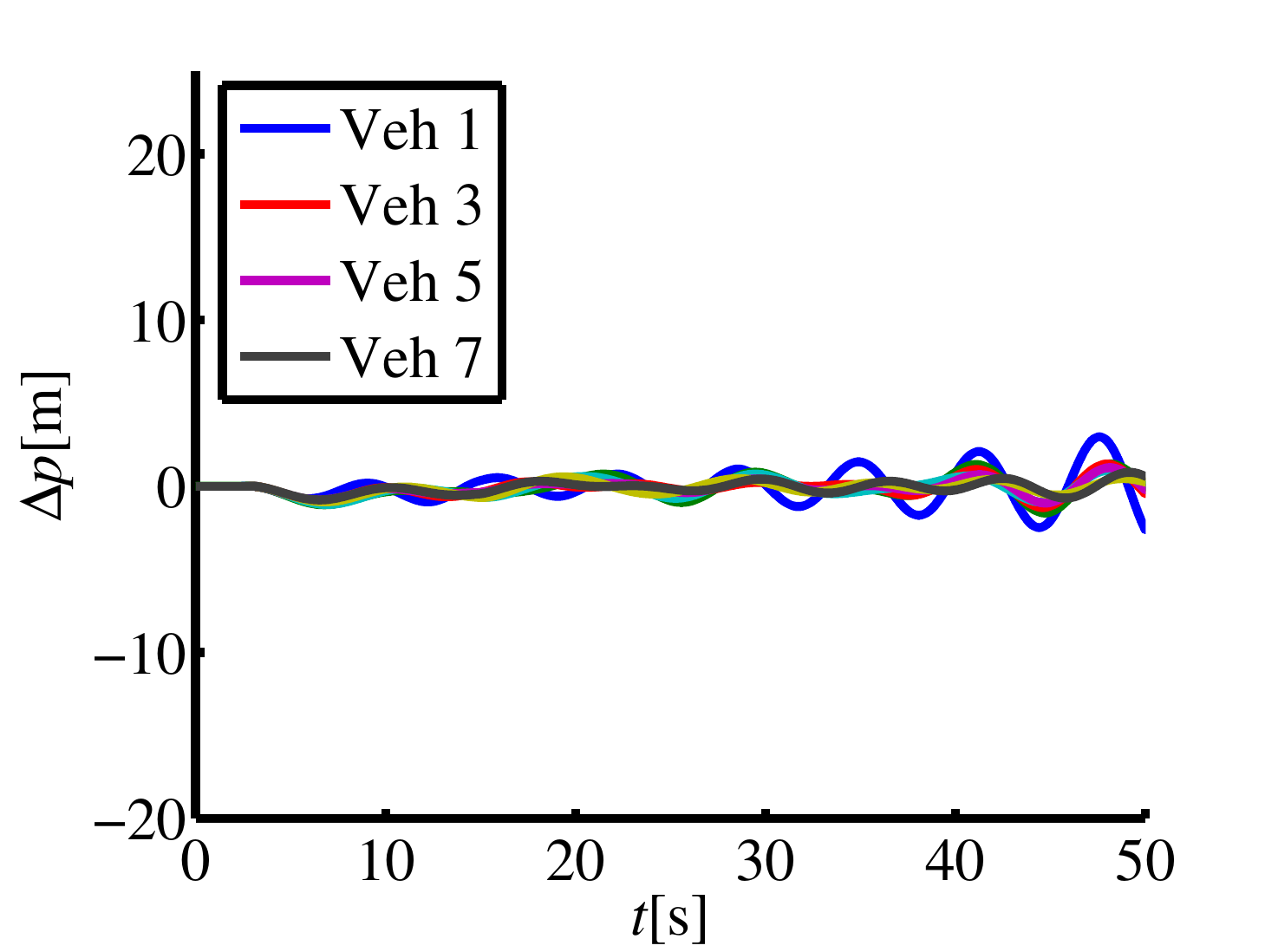}}
  \caption{Spacing error profiles using the feedback gains $(k_{ip},\hat{k}_{iv},k_{ia})$ in Table~\ref{tab:1} (Theorem \ref{theorem:2} is not satisfied) and the velocity profile \eqref{eq:39} for heterogeneous platoons with different topologies: (a) PF, (b) PLF, (c) TPF, (d) TPLF.}
  \label{fig:5}
\end{figure}

First, we validate the asymptotical stability given in Theorem~\ref{theorem:2}. Given the heterogeneous inertial time lags $\tau_i$ in Table~\ref{tab:1}, it is easy to check that the feedback gains $(k_{ip},k_{iv},k_{ia})$ satisfy the conditions in Theorem~\ref{theorem:2} for all the four types of DAGs shown in Fig.~\ref{fig:3}, while the feedback gains $(k_{ip},\hat{k}_{iv},k_{ia})$ do not. The spacing error profiles\footnote{Throughout the paper, we plotted the spacing error profiles of all vehicles in the figures, but only the legends of odd numbered vehicles were given due to space limit.} corresponding to the two sets of feedback gains are shown in Fig.~\ref{fig:4} and Fig.~\ref{fig:5}, respectively. The numerical results clearly demonstrate that the feedback gains $(k_{ip},k_{iv},k_{ia})$ can guarantee the asymptotical stability, while the feedback gains $(k_{ip},\hat{k}_{iv},k_{ia})$ cannot. This fact supports the statements in Theorem~\ref{theorem:2}. Let us denote the maximum spacing error as $\Delta p_{\max}$ during the transient process. As shown in Fig.~\ref{fig:4}, we have observed that $\Delta p_{\max,\text{PF}}>\Delta p_{\max,\text{TPF}}>\Delta p_{\max,\text{PLF}}\approx\Delta p_{\max,\text{TPLF}}$ for vehicles $2\sim7$. Moreover, similar patterns can be found in the subsequent simulations. The simulation results suggest that 1) more predecessor's information has the potential to improve the tracking performance of platoons, which is consistent with the concept of \emph{look ahead control}~\cite{turri2017cooperative}; 2) the leader's information plays a key role in regulating local behavior of vehicles in a platoon, which agrees with the theoretical analysis on homogeneous platoons~\cite{zheng2017platooning,zheng2015stabilityMargin}.

\begin{figure}[t]
  \centering
  \setlength{\abovecaptionskip}{0pt}
  \setlength{\belowcaptionskip}{0em}
  \subfigure[ ]{\includegraphics[width=0.45\columnwidth]{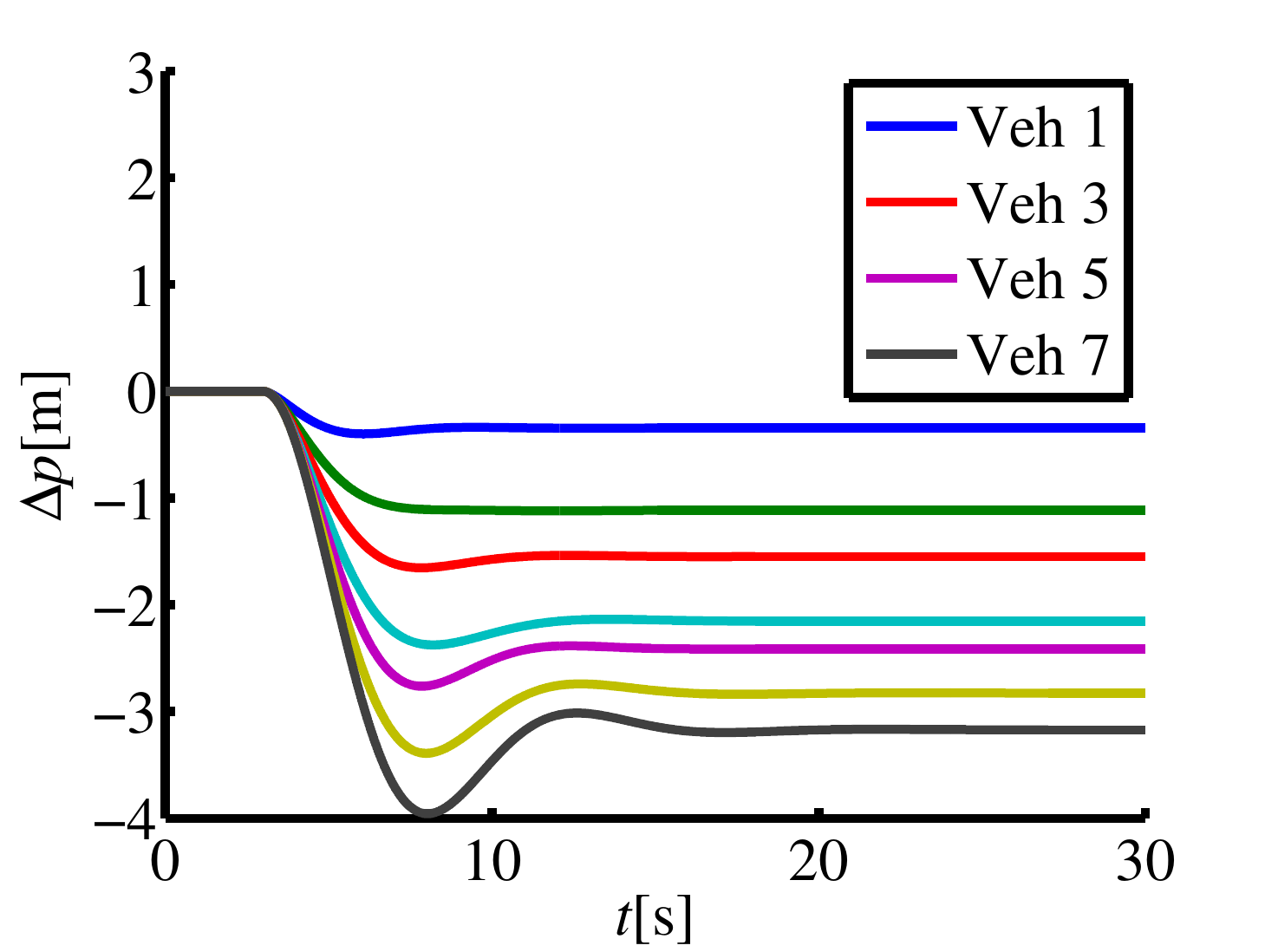}}
  \subfigure[ ]{\includegraphics[width=0.45\columnwidth]{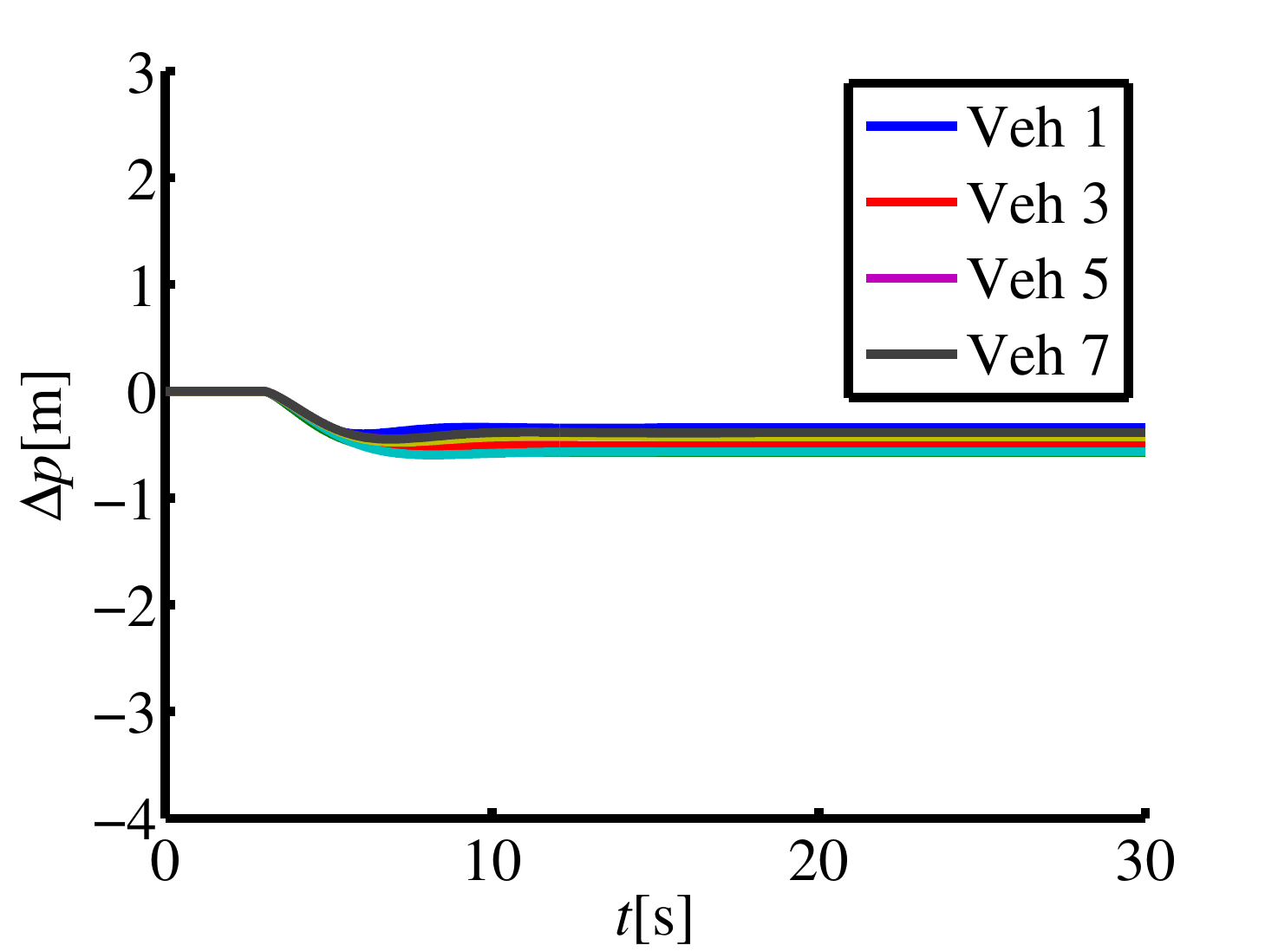}}
  \subfigure[ ]{\includegraphics[width=0.45\columnwidth]{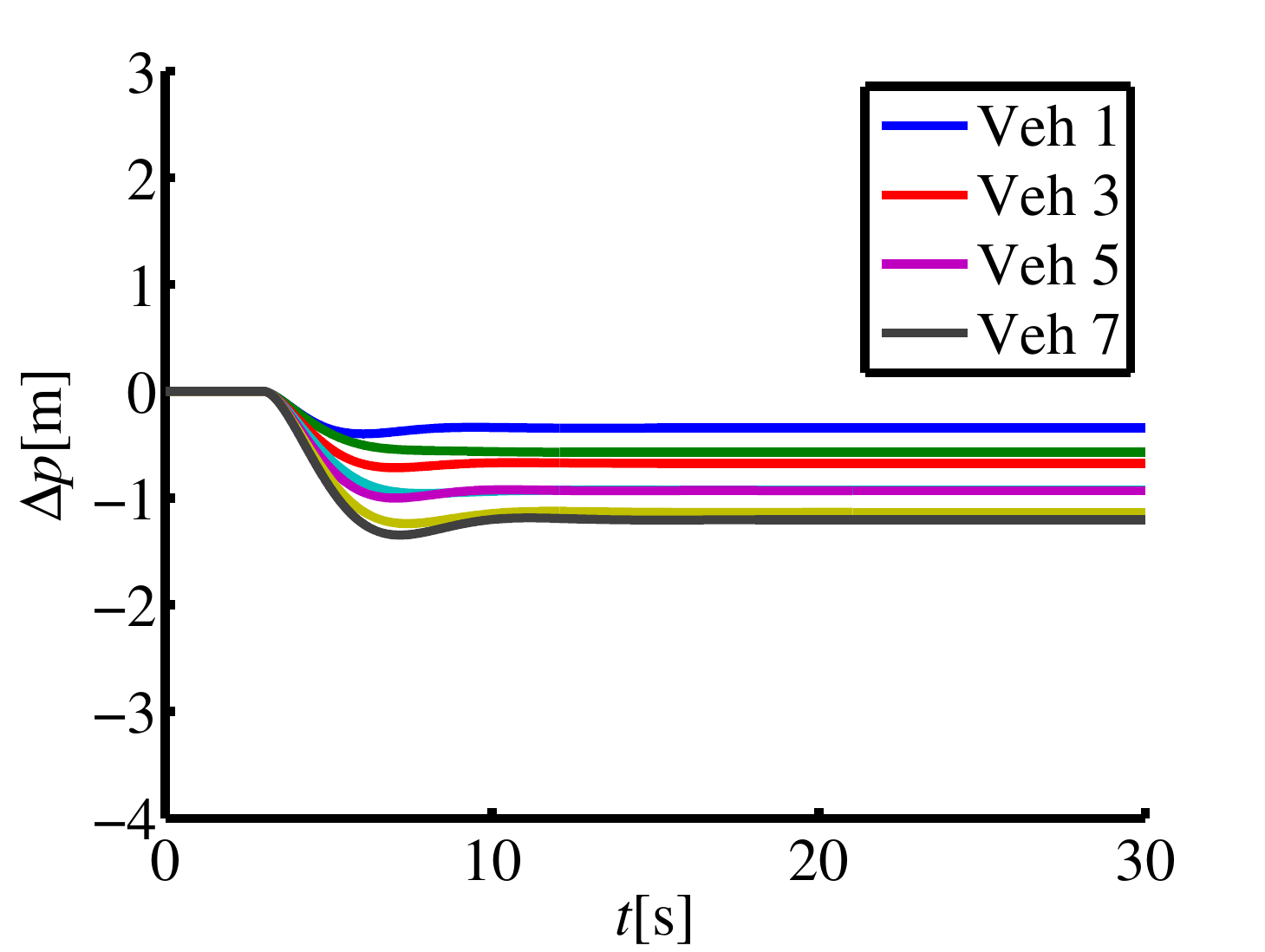}}
  \subfigure[ ]{\includegraphics[width=0.45\columnwidth]{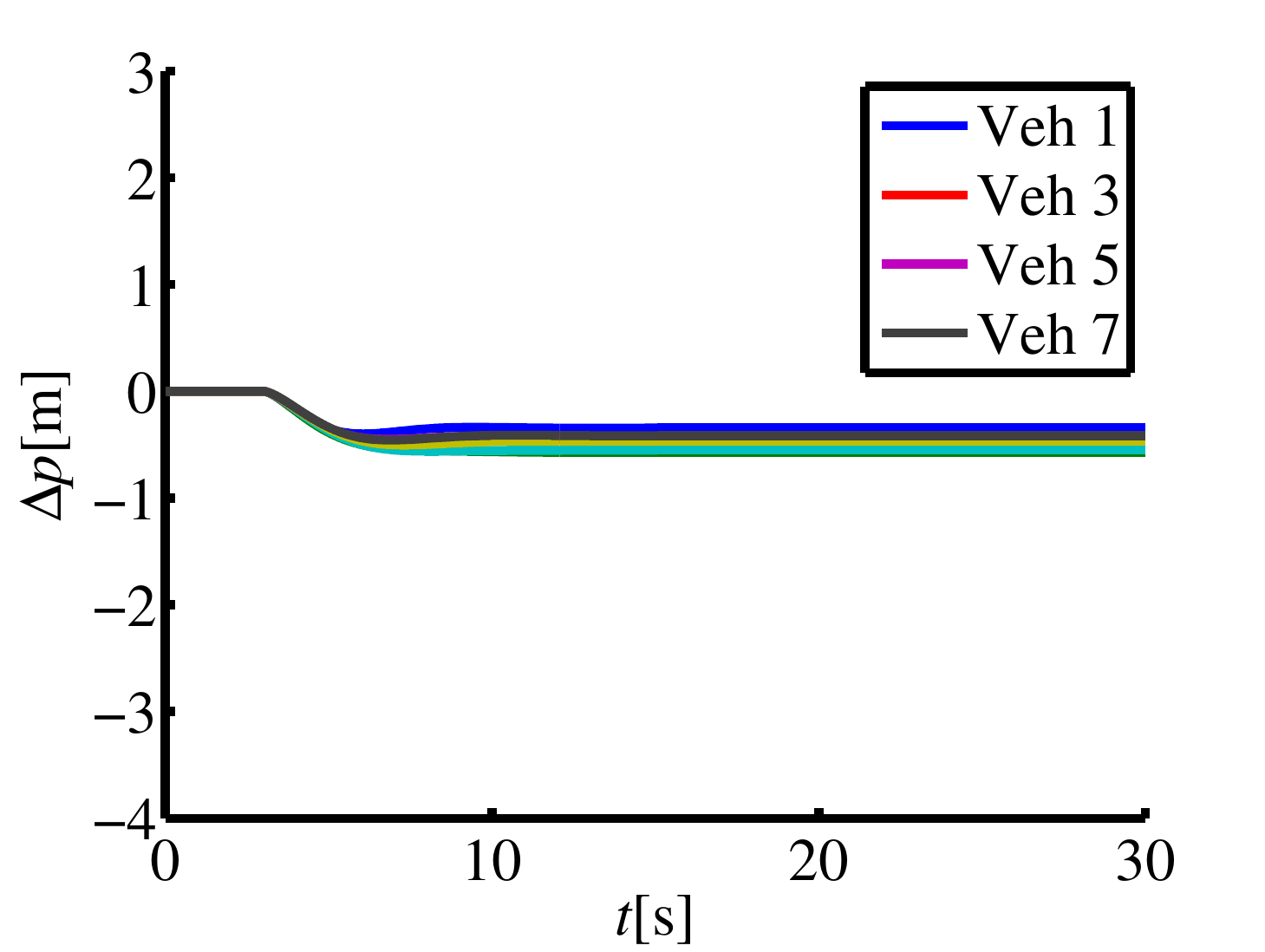}}
  \caption{Spacing error profiles using feedback gains $(k_{ip},k_{iv},k_{ia})$ (Theorem \ref{theorem:2} is satisfied) and velocity profile \eqref{eq:40} for heterogeneous platoons with different topologies: (a) PF, (b) PLF, (c) TPF, (d) TPLF.}
  \label{fig:6}
\end{figure}

Second, we show that it is impossible to track a leading vehicle with acceleration (see Theorem~\ref{theorem:1}). Here, we extend the acceleration process of $v_0(t)$ in Eq.~\eqref{eq:39}, and consider a scenario where the speed profile of the leader is
\begin{equation}\label{eq:40}
    v_0(t)=\begin{cases}
            10 & 0 s\leq t<3 s \\
            10+t & t\geq 3 s\\
        \end{cases},(m/s).
\end{equation}
The spacing error profiles using control gains $(k_{ip},k_{iv},k_{ia})$ and velocity profile~\eqref{eq:40} are shown in Fig. \ref{fig:6}. It is clear that there exists a constant spacing error for each following vehicle in the platoon, meaning that the followers are not able to track the leader if the leader continues to accelerate. When the leader accelerates from one speed to another speed and maintains that speed afterwards, the followers would be able to reach consensus to the new steady state defined by the leader. The transient behavior of each follower is normally studied under the notion of string stability~\cite{ploeg2014controller}, which is beyond the scope of the current manuscript.

\begin{table}[t]
    \centering
    \renewcommand\arraystretch{1.2}
    \caption{Performance index values corresponding to different controller parameters and information flow topologies}
    \begin{tabular*}{0.8\columnwidth}{c @{\extracolsep{\fill}}  c c c c }
        \hline \toprule[1pt]
        \multirow{2}{*}{$\varepsilon_i$} &  \multicolumn{4}{c}{$T_c[s]$} \\
        \cline{2-5}
            & PF & PLF & TPF & TPLF \\
        \hline
        1 & 23.71 & 18.27 & 18.71 & 18.29 \\
        3 & 21.89 & 17.42 & 18.14 & 17.44 \\
        5 & 20.94 & 17.07 & 17.90 & 17.09 \\
        7 & 19.95 & 16.85 & 17.73 & 16.87 \\
        \bottomrule[1pt]
    \end{tabular*}
    \label{tab:2}
\end{table}

\begin{figure}[t]
  \centering
  \setlength{\abovecaptionskip}{0pt}
  \setlength{\belowcaptionskip}{0em}
  \subfigure[ ]{\includegraphics[width=0.45\columnwidth]{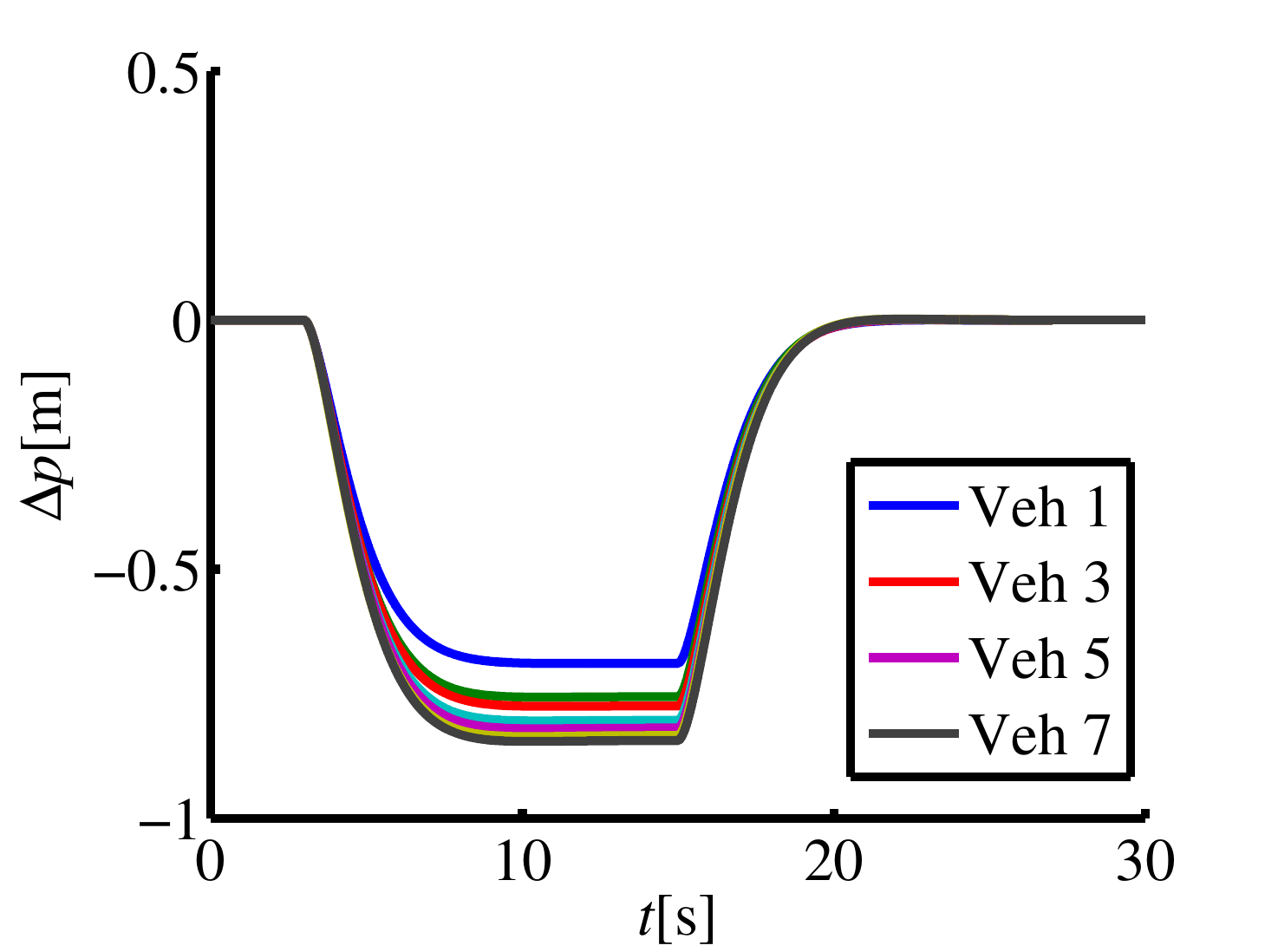}}
  \subfigure[ ]{\includegraphics[width=0.45\columnwidth]{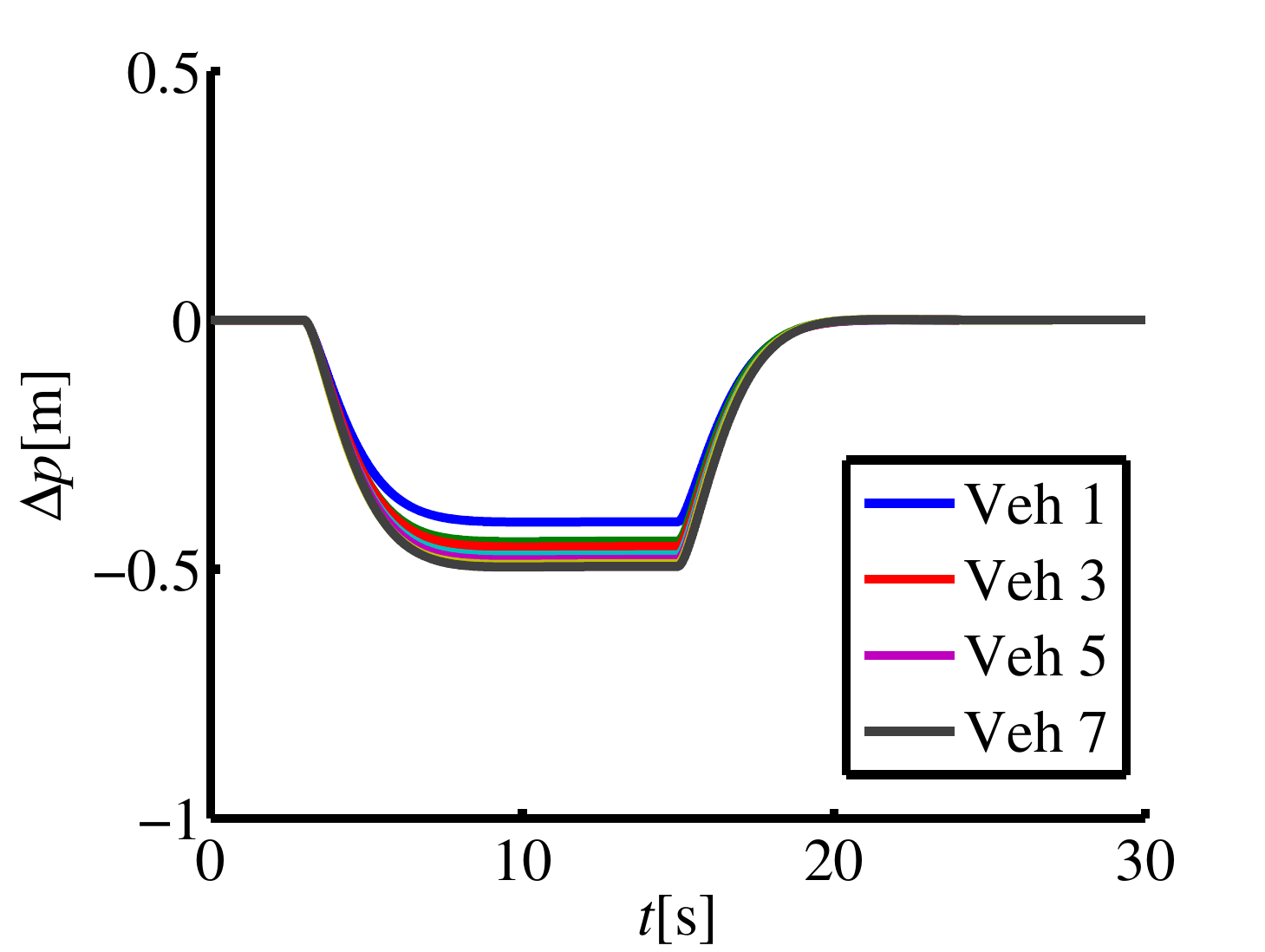}}
  \subfigure[ ]{\includegraphics[width=0.45\columnwidth]{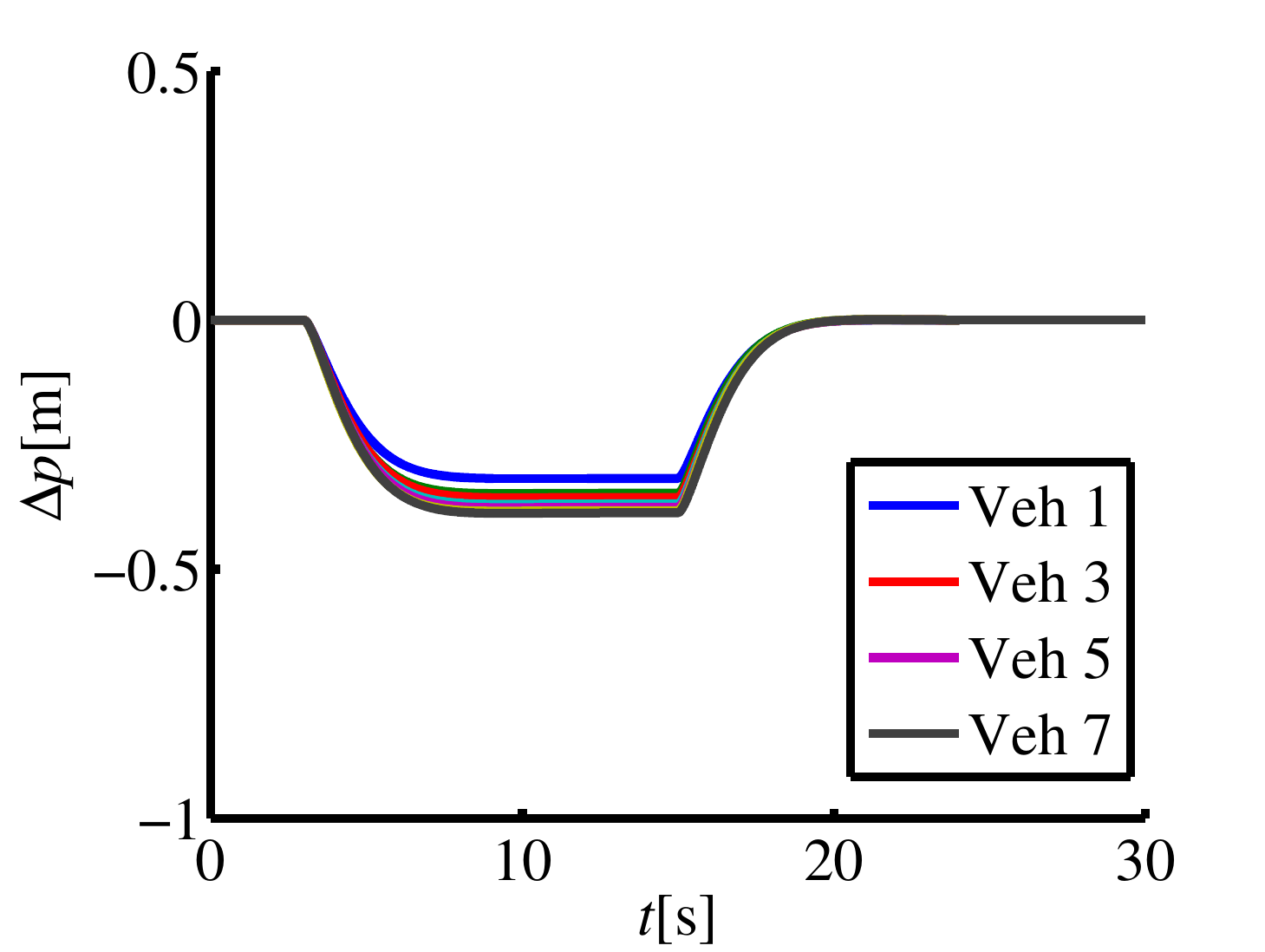}}
  \subfigure[ ]{\includegraphics[width=0.45\columnwidth]{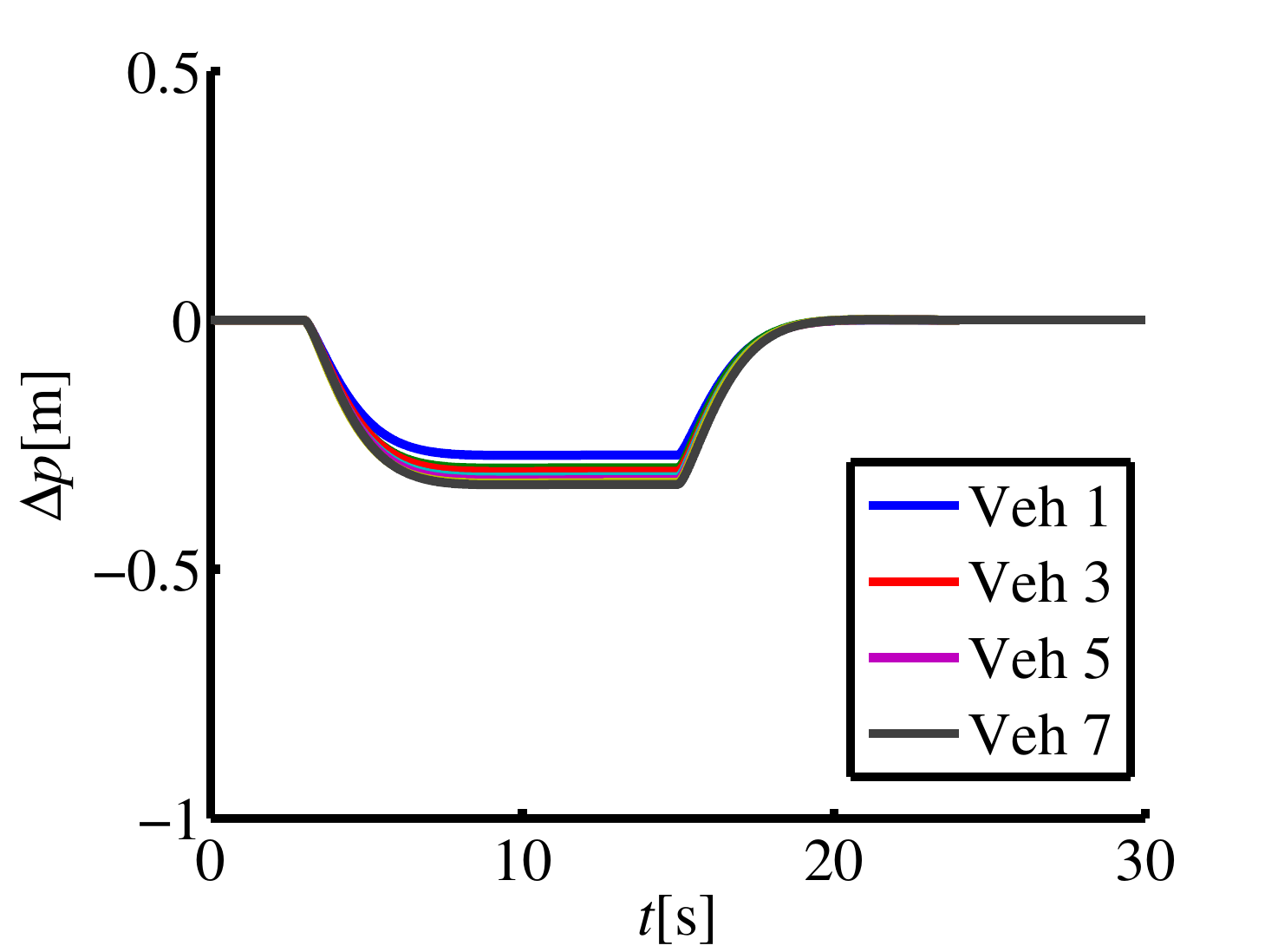}}
  \caption{Spacing error profiles for platoons with TPLF topology using the design method of feedback gains given in Theorem~\ref{theorem:3}: (a) $\varepsilon_i=1$, (b) $\varepsilon_i=3$, (c) $\varepsilon_i=5$, (d) $\varepsilon_i=7$.}
  \label{fig:7}
\end{figure}

Third, we demonstrate that the stability is guaranteed using the design method in Theorem \ref{theorem:3} and the convergence rate can be adjusted by tuning the parameter $\varepsilon_i$. Here we define an index to quantify the convergence rate, as used in~\cite{zheng2015stabilityMargin},
\begin{equation}\label{eq:41}
  T_c=\min_{T_1}\left(\max_{i\in\mathcal{N},t>T_1}|\hat p_i(t)|<\delta\right),(s)
\end{equation}
where $\delta$ is a threshold for admissible position tracking error. This performance index $T_c$ indicates the time instant when all the following vehicles reach consensus with the leading vehicle with respect to an admissible error $\delta$. In the simulations, we set $\delta=0.1 m$ and $$\displaystyle \alpha_i=\frac{1}{2(d_{ii}+p_{ii})}+1,$$ and we used the velocity profile~\eqref{eq:39}. The performance index values corresponding to different $\varepsilon_i$ (we used identical $\varepsilon_i$ for each vehicle) and information topologies are listed in Table \ref{tab:2}. The spacing error profiles for TPLF topology are given in Fig.~\ref{fig:7}, which shows that the spacing error converges to $0$, thus the closed-loop system is stable. In Table \ref{tab:2}, it is obvious that the convergence rate can be improved by increasing $\varepsilon_i$. A trade-off is that a larger $\varepsilon_i$ generally leads to larger feedback gains, which might result in actuator saturations. In practice, it is necessary to tune $\varepsilon_i$ to obtain a suitable feedback controller for a particular platoon system.

\subsection{Simulations with a nonlinear model}

\begin{table}[t]
    \centering
    \renewcommand\arraystretch{1.2}
    \caption{Model parameters of the nonlinear heterogeneous platoon}
    \begin{tabular*}{0.8\columnwidth}{c c c }
        \hline \toprule[1pt]
        parameter & true value & estimated value \\
        \hline
        $m_i$ [kg]      & 1500 + 100$\times{i}$  & 1700\\
        $\tau_i$    & 0.30 + 0.02$\times  i$      & 0.34\\
        $\eta_i$    & 0.80 + 0.01$\times  i$      & 0.82\\
        $C_{A,i}$ [kg/m]  & 0.40 + 0.01$\times i$ & 0.42\\
        $r_i$ [m]      & 0.250 + 0.005$\times i$  & 0.26\\
        $f_i$       & 0.015 + 0.001$\times i$    & 0.017\\
        \bottomrule[1pt]
    \end{tabular*}
    \label{tab:3}
\end{table}

As our final experiment, we validate the effectiveness of the proposed control method using a realistic nonlinear vehicle model. As used in~\cite{Zheng2016distributed}, we consider the following nonlinear dynamics:
\begin{eqnarray}\label{eq:43}
  \begin{cases}
  \dot{p}_i(t)=v_i(t), \\
  m_i\dot{v}_i(t)=\frac{\eta_{i}}{r_i}T_i(t)-C_{A,i}v^2_i(t)-m_igf_i,\\
  \tau_i\dot{T}_i(t)+T_i(t)=T_{\text{des},i}(t),
  \end{cases}\forall i\in\mathcal{N},
\end{eqnarray}
where $p_i$ and $v_i$ are the position and velocity; $T_{\text{des},i}$, $T_i$, and $\eta_{i}$ are the desired driving torque (control input), the actual driving torque, and the mechanical efficiency of the driveline; $m_i$ and $r_i$ are the mass and tire radius; $C_{A,i}$, $f_i$, and $g$ are the coefficients of the aerodynamics drag, rolling resistance, and gravitational acceleration.

The true and estimated values of the model parameters are listed in Table~\ref{tab:3}. To improve the robustness to the parameter uncertainty, we employ the integral sliding mode control~\cite{guldner2009sliding} by adding a non-smooth term to the nominal control input $u_i$. Thus, the new robust control input $\tilde{u}_i$ and the feedback linearization law are
\begin{equation} \label{eq:44}
    \begin{aligned}
    \tilde{u}_i &=u_i-k_{s,i}\text{sign}(s_i), \\
    T_{\text{des},i}&=\frac{\hat{r}_i}{\hat{\eta}_i}(\hat{m}_{i}\hat{u}_i+\hat m_i\hat f_ig+2\hat C_{A,i}\hat \tau_iv_ia_i+\hat C_{A,i}v_i^2),
    \end{aligned}
\end{equation}
where $k_{s,i}$ is the feedback gain, $\text{sign}(\cdot)$ is the sign function, $\hat\theta_i$ is the estimation of the model parameter $\theta_i$, $s_i$ is the integral sliding mode variable, of which the derivative is
\begin{equation*}\label{eq:45}
\dot{s}_i=\hat{\tau}_{i}\dot{a}_{i}+a_{i}-u_{i},
\end{equation*}
and the integral initial value is zero.

\begin{figure}[t]
  \centering
  \setlength{\abovecaptionskip}{0pt}
  \setlength{\belowcaptionskip}{0em}
  \subfigure[ ]{\includegraphics[width=0.45\columnwidth]{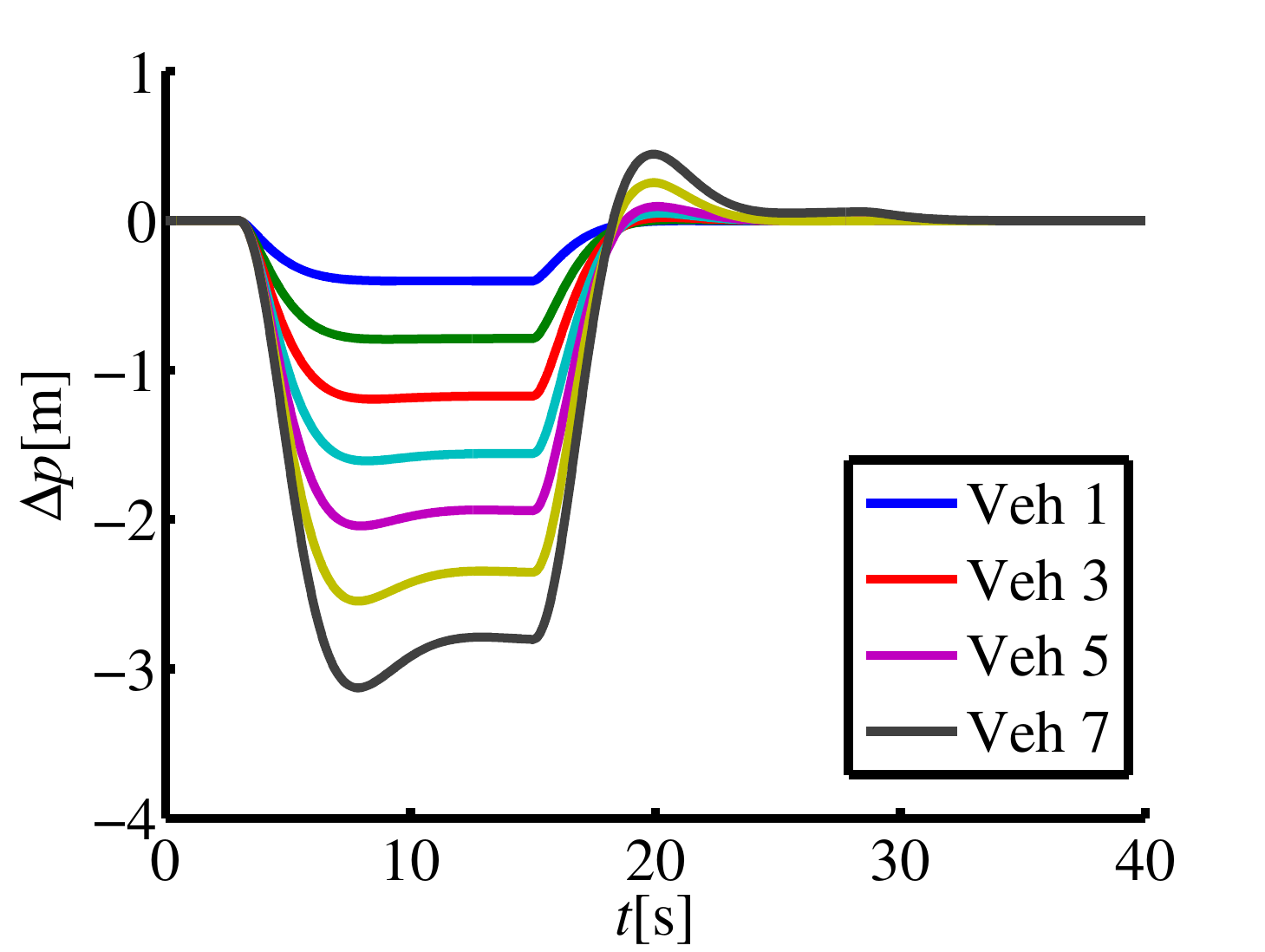}}
  \subfigure[ ]{\includegraphics[width=0.45\columnwidth]{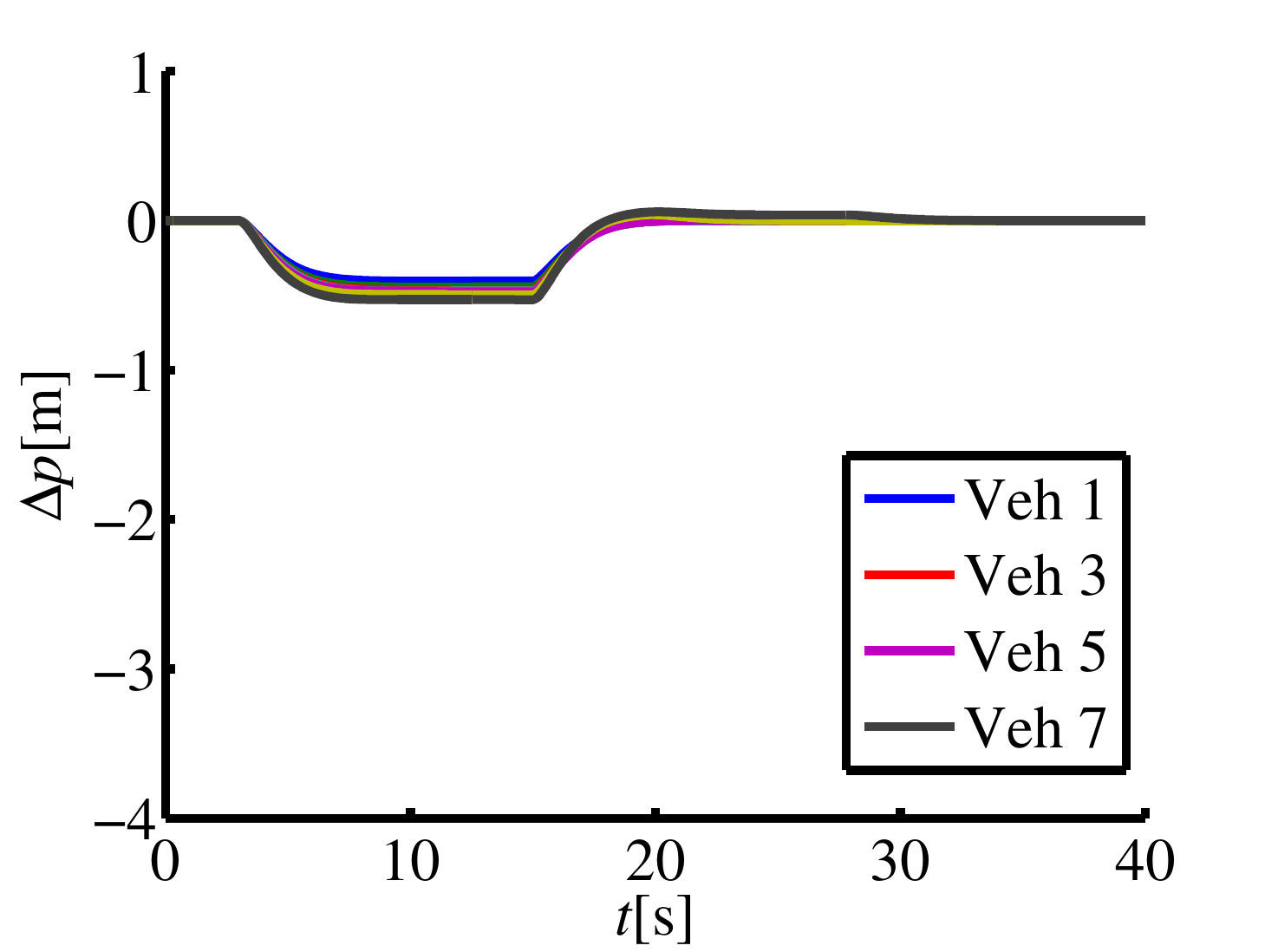}}
  \subfigure[ ]{\includegraphics[width=0.45\columnwidth]{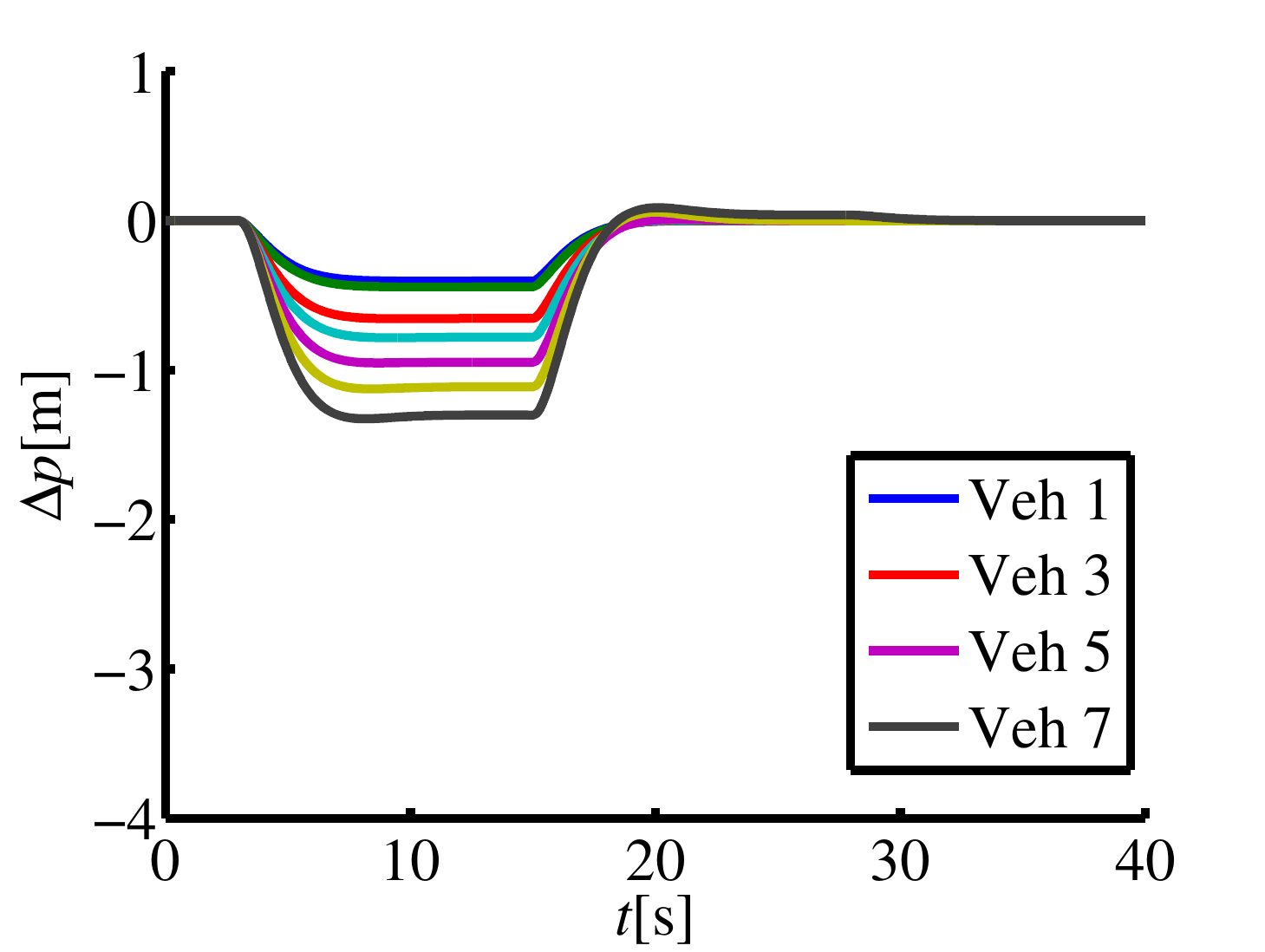}}
  \subfigure[ ]{\includegraphics[width=0.45\columnwidth]{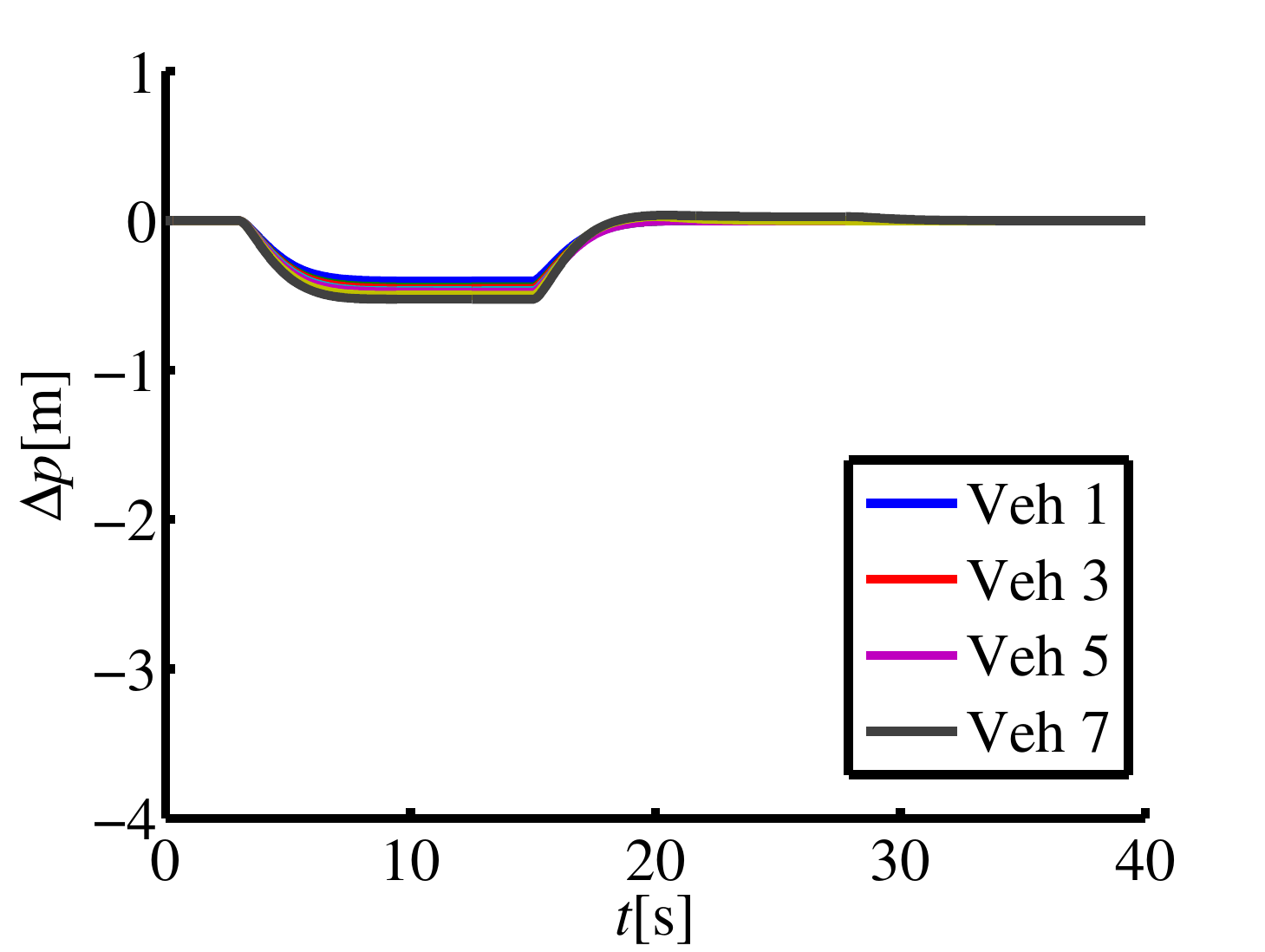}}
  \caption{Spacing error profiles using the robust controller and velocity profile \eqref{eq:39} for nonlinear heterogeneous platoons with different topologies: (a) PF, (b) PLF, (c) TPF, (d) TPLF.}
  \label{fig:8}
\end{figure}

Here we use the method in Theorem \ref{theorem:3} to synthesize the nominal controller $u_i$ by taking $\varepsilon_i=3,\forall i\in\mathcal{N}$. For the design of the robust controller $\hat{u}_i$, we take $k_{s,i}=0.3, \forall i\in\mathcal{N}$. With the velocity profile~\eqref{eq:39}, the spacing error profiles for different topologies are given in Fig. \ref{fig:8}. It is clear that the ARE-based robust controller can guarantee the asymptotic stability for this nonlinear heterogeneous platoon.

\section{Conclusion}\label{section:6}

In this paper, we have studied the cooperative control of heterogeneous connected vehicles with directed acyclic interactions. The problem formulation has directly considered heterogeneity in both dynamics and feedback design. We applied the internal model principle and exploited the lower-triangular structure in the closed-loop system. The main conclusions are: (1) we analyzed the tracking ability using the internal model principle, and proved that the leader should run at a constant speed in order to reach a consensus state (\emph{i.e.}, Theorem~\ref{theorem:1}); (2) under the notion of directed acyclic interactions, we derived a set of necessary and sufficient conditions to stabilize heterogeneous platoons (\emph{i.e.}, Theorem~\ref{theorem:2}); (3) we introduced an algebraic Riccati equation based method for synthesizing feedback gains, which has a relatively clear physical interpretation (\emph{i.e.}, Theorem~\ref{theorem:3}). The results offered some insights on the influence of heterogeneity on the collective behavior of connected vehicles.

One future work is to consider the influence of different inter-vehicle spacing policies on the stability region of feedback gains. String stability of heterogeneous platoons with various communication topologies is another important topic that is worth further investigation. Future work will also address the effect of imperfect communication, such as time delay and packet loss, on the platoon control. In addition, it would be interesting to investigate other decomposition methods, such as chordal decomposition~\cite{zheng2018scalable,vandenberghe2015chordal}, in the control of multiple connected vehicles, which may offer another effective way to deal with the heterogeneity issue.
\balance
\bibliographystyle{IEEEtran}
\bibliography{Reference}

\end{document}